\documentclass[a4paper,12pt,twoside]{article}
\usepackage{graphicx}\usepackage{a4wide}
\usepackage[pdftex,colorlinks,bookmarksnumbered,bookmarksopen,citecolor=red,urlcolor=red]{hyperref}
\usepackage{amssymb,amsmath,amsthm}
\usepackage{subfigure,color,colordvi,graphicx,xcolor}
\usepackage[only,llbracket,rrbracket,interleave]{stmaryrd}
\usepackage[sort,compress]{cite}
\usepackage{multirow}

\newcommand{\bm}[1]{\text{\boldmath $#1$\unboldmath}}
\newcommand{\vect}[1]{\mathbf{#1}}
\newcommand{\mat}[1]{\mathbf{#1}}
\newcommand{\grad}{\bm{\nabla}}
\newcommand{\RR}{\mathbb{R}}

\newcommand{\VhHat}{\ensuremath{\mathcal{\hat{V}}^h}}
\newcommand{\Vh}{\ensuremath{\mathcal{V}^h}}

\newcommand{\eltwo}{\ensuremath{\mathcal{L}_2}}

\newcommand{\ndof}  {\ensuremath{\texttt{n}_{\texttt{dof}}}}

\newcommand{\nsd}  {\ensuremath{\texttt{n}_{\texttt{sd}}}}
\newcommand{\msd}  {\ensuremath{\texttt{m}_{\texttt{sd}}}}
\newcommand{\numel}{\ensuremath{\texttt{n}_{\texttt{el}}}}
\newcommand{\numfa}{\ensuremath{\texttt{n}_{\texttt{fa}}}}

\newcommand{\hv}{\hat{v}}

\newcommand{\bu}{\bm{u}}
\renewcommand{\bf}{\bm{f}}

\newcommand{\bhw}{\widehat{\bw}}
\newcommand{\bhu}{\widehat{\bu}}

\newcommand{\bv}{\bm{v}}

\newcommand{\bn}{\bm{n}}
\newcommand{\bx}{\bm{x}}
\newcommand{\bL}{\bm{L}}

\newcommand{\btau}{\bm{\tau}}

\newcommand{\Le}{\mat{L}_e}
\newcommand{\bue}{\vect{u}_e}

\newcommand{\buHi}{\hat{\vect{u}}_i}
\newcommand{\buHj}{\hat{\vect{u}}_j}

\newcommand{\nDeg}{\ensuremath{k}}
\newcommand{\Pk}{\ensuremath{\mathcal{P}^{\nDeg}}}

\newcommand{\Insd}{\mat{I}_{\nsd}}

\newcommand{\stress}{\bm{\sigma}}
\newcommand{\defo}[1]{\bm{\varepsilon}(#1)}
\newcommand{\gradS}{\bm{\nabla}_{\texttt{S}}}
\newcommand{\stressV}{\bm{\sigma}_{\texttt{V}}}
\newcommand{\strainV}{\bm{\varepsilon}_{\texttt{V}}}
\newcommand{\bD}{\mat{D}}

\newcommand{\bV}{\mat{V}}
\newcommand{\bLambda}{\mat{\Lambda}}
\newcommand{\bDHalf}{\widetilde{\bD}}
\newcommand{\bN}{\mat{N}}

\newcommand{\bg}{\bm{g}}
\newcommand{\elasTensor}{\mat{C}}
\newcommand{\Pn}{\mat{P}_n}
\newcommand{\Pt}{\mat{P}_t}

\newcommand{\jump}[1]{\llbracket #1\rrbracket}

\newcommand{\bw}{\bm{w}}

\newcommand{\Aset}{\mathcal{A}_e}
\newcommand{\Bset}{\mathcal{B}_e}
\newcommand{\Dset}{\mathcal{D}_e}
\newcommand{\Nset}{\mathcal{N}_e}

\newcommand{\Iset}{\mathcal{I}_e}

\newenvironment{keywords}{\begin{quote}\emph{\textbf{Keywords:}}}{\end{quote}}

\theoremstyle{definition}

\newtheorem{remark}{Remark}

\graphicspath{{figs/}}

\begin{document}
%==========================================================================
\title{A locking-free face-centred finite volume (FCFV) method for linear elasticity}

\author{Ruben Sevilla\\[-1ex]
             \small Zienkiewicz Centre for Computational Engineering, \\[-1ex]
             \small College of Engineering, Swansea University, Wales, UK \\[1em]
             Matteo Giacomini, and Antonio Huerta\\[-1ex]
             \small Laboratori de C\`alcul Num\`eric (LaC\`aN), \\[-1ex]
             \small ETS de Ingenieros de Caminos, Canales y Puertos, \\[-1ex]
             \small Universitat Polit\`ecnica de Catalunya, Barcelona, Spain}
\date{\today}
%________________________________________________________________________
\maketitle

%==========================================================================
\begin{abstract}
A face-centred finite volume (FCFV) method is proposed for the linear elasticity equation. The FCFV is a mixed hybrid formulation, featuring a system of first-order equations, that defines the unknowns on the faces (edges in two dimensions) of the mesh elements. The symmetry of the stress tensor is strongly enforced using the well-known Voigt notation and the displacement and stress fields inside each cell are obtained element-wise by means of explicit formulas. 
The resulting FCFV method is robust and locking-free in the nearly incompressible limit. Numerical experiments in two and three dimensions show optimal convergence of the displacement and the stress fields without any reconstruction. Moreover, the accuracy of the FCFV method is not sensitive to mesh distortion and stretching. Classical benchmark tests including Kirch's plate and Cook's membrane problems in two dimensions as well as three dimensional problems involving shear phenomenons, pressurised thin shells and complex geometries are presented to show the capability and potential of the proposed methodology.
\end{abstract}

%________________________________________________________________________

%==========================================================================
\section{Introduction}
\label{sc:Intro}
%==========================================================================

Despite the finite volume method (FVM) was originally proposed in the context of hyperbolic systems of conservation laws~\cite{McDonald71,RizziInouye73}, there has been a growing interest towards its application to other physical problems, including the simulation of deformable structures~\cite{WHEEL1996311,FAINBERG1996167}. Several robust and efficient implementations of the FVM are available in both open-source and commercial libraries, making it an extremely attractive approach for industrialists.

The existing finite volume paradigms discussed in the structural mechanics community can be classified into two families, depending on the localisation of the unknowns in the computational mesh: the cell-centred finite volume (CCFV) method~\cite{doi:10.1002/nme.1620372110,bijelonja2006finite,lee2013development,Haider-HLGHB:18} defines the unknowns at the centroid of the mesh elements, whereas the vertex-centred finite volume (VCFV) strategy~\cite{SLONE200369,Xia2006,SULIMAN20142265} sets the unknowns at the mesh nodes.
A major limitation of the CCFV method is the poor approximation of the gradient of the displacements at the faces using unstructured meshes~\cite{doi:10.1002/nme.1620381010,cardiff2016block}. To overcome this issue, Jasak and Weller~\cite{JasakWeller2000} proposed a correction to match the value of the gradient of the displacements on the faces to the neighbouring cells. More recently, Nordbotten and co-workers enriched the classical CCFV formulation with a discrete expression of the stresses on the mesh faces~\cite{doi:10.1002/nme.4734,doi:10.1137/140972792} and the resulting multi-point stress approximation was shown to improve the description of the stresses at the interface between two cells.
Similar drawbacks are experienced by the VCFV strategy which also requires a reconstruction of the gradient of the displacements to guarantee the first-order convergence of the stresses.
Within this context, the accuracy of the reconstruction may suffer from the non-orthogonality of the mesh and poor approximations may result from the use of highly deformed grids.

Another critical aspect in the numerical treatment of linear elasticity problems is represented by the fulfilment of the balance of angular momentum which implies the symmetry of the stress tensor~\cite{MR2449101}. Starting from the pioneering work of Fraejis de Veubeke~\cite{Veubeke1975}, finite element formulations with weakly enforced symmetry of the stress tensor have been extensively studied in the literature~\cite{brezzi1991mixed}. In~\cite{keilegavlen2017finite}, the weak imposition of the symmetry of the stress tensor is investigated in the context of a cell-centred finite volume paradigm.

Recently, a face-centred finite volume (FCFV) method which defines the unknowns over the faces of the mesh elements has been introduced for Poisson and Stokes problems~\cite{doi:10.1002/nme.5833}. In the present work, the FCFV method is extended to simulate the behaviour of deformable bodies under the assumption of small displacements. Starting from the hybridisable discontinuous Galerkin (HDG) method by Cockburn and co-workers~\cite{Jay-CGL:09,Cockburn-CDG:08,soon2009hybridizable,MR3340089}, the discrete finite volume system is derived by setting a constant degree of approximation in the recently proposed HDG formulation of the linear elasticity equation based on Voigt notation~\cite{HDGVoigtElasticity}.
The resulting finite volume strategy involves the solution of a symmetric system of equations to determine the displacements on the mesh faces (edges in two dimensions). The displacement and the stress fields inside each element are then retrieved via explicit closed expressions defined element-by-element. The enforcement of the symmetry of the stress tensor via the Voigt notation allows to strongly fulfil the balance of angular momentum and to obtain optimal convergence for both the displacement and stress fields without any reconstruction. Therefore, the solution of the FCFV method does not deteriorate in presence of highly stretched or distorted elements. In addition, it is worth emphasising that other HDG methods reported in the literature (e.g. \cite{soon2009hybridizable,MR3340089}), without the Voigt notation proposed in this paper, have shown a sub-optimal rate of convergence in the stress field.

Special attention is given to elastic problems in which classical numerical methods experience volumetric or shear locking. Locking-free finite volume formulations for bending plates~\cite{cook2001concepts} have been discussed for both cell-centred and vertex-centred formulations by Wheel~\cite{WHEEL1997199} and Fallah~\cite{FALLAH20043457}. Nevertheless, using solid elements, VCFV approaches experience shear locking and additional rotational degrees of freedom are required to handle rigid body motions and accurately predict membrane deformations~\cite{wenke2003finite,PAN20101506}. In the nearly incompressible limit, the proposed FCFV method is locking-free and the optimal convergence properties are preserved for both the displacement and the stress fields.

The remaining of this paper is organised as follows. In Section~\ref{sc:problem}, the linear elasticity equation using Voigt notation is briefly recalled. The proposed FCFV scheme is presented in Section~\ref{sc:FCFV}. Section~\ref{sc:examples2D} is devoted to the numerical validation of the method in two dimensions. In particular, the optimal orders of convergence are checked for the displacement and stress fields, a sensitivity analysis to the stabilisation parameter and the mesh distortion is performed and the locking-free behaviour is verified for nearly incompressible materials using Kirch's plate and Cook's membrane test cases. In Section~\ref{sc:examples32D}, several three-dimensional problems involving shear phenomenons, pressurised thin shells and complex geometries under realistic loads are discussed to show the capability of the method to handle complex geometries. Finally, Section~\ref{sc:Conclusion} summarises the conclusions of the work that has been presented.

%==========================================================================
\section{Problem statement}
\label{sc:problem}
%==========================================================================

Given an open bounded domain $\Omega \subset \mathbb{R}^{\nsd}$, where $\nsd$ denotes the number of spatial dimensions, the boundary $\partial\Omega$ is partitioned into the non-overlapping Dirichlet and generalised Neumann boundaries, $\Gamma_D$ and $\Gamma_N$ respectively. The behaviour of a deformable solid medium $\Omega$ is described by
\begin{equation} \label{eq:elasticity}
\left\{\begin{aligned}
-\grad \cdot \stress  &= \bf       &&\text{in $\Omega$,}\\
\bu &= \bu_D  &&\text{on $\Gamma_D$,}\\
(1 - \xi) \Pn \bu +  (\Pt + \xi \Pn) \bn\cdot\stress &  = \xi \bg, &&\text{on $\Gamma_N$,}\\
\end{aligned}\right.
\end{equation}
where $\stress$ is the Cauchy stress tensor, $\bf$ is the external force, $\bu$ is the displacement field vector, $\bn$ is the outward unit normal vector to $\Gamma_N$ and the normal and tangent projection matrices are defined as $\Pn = \bn \otimes \bn$ and $\Pt = \Insd - \bn \otimes \bn$ respectively. The boundary conditions are given by the imposed displacements on the Dirichlet boundary, $\bu_D$, and the traction vector on the Neumann boundary, $\bg$. The parameter $\xi$ can take a value of one for a pure Neumann boundary or zero for an artificial symmetry boundary, where the normal displacement and the tangential tractions vanish. 

\begin{remark}
It is worth noting that a more general boundary condition can be considered to include Dirichlet, Neumann and symmetry boundaries. However, due to the different treatment of Dirichlet boundary conditions in the proposed numerical methodology, the form stated in Equation~\eqref{eq:elasticity} is preferred in this work.
\end{remark}

For a linear elastic material, the well-known Hooke's law provides the relation between stress and strain, namely $\stress = \elasTensor\!:\!\defo{\bu}$, where $\elasTensor$ is the fourth order elasticity tensor and the linearised strain tensor is $\defo{\bu} := \left( \grad \bu + \grad \bu^T \right)/2$. 

The so-called Voigt notation~\cite{FishBelytschko2007} is common in this context. The main idea is to exploit the symmetry of the strain and stress tensors. To this end, the strain and stress tensors are reduced to vectors by storing only the non-redundant terms, namely $\strainV := \left[\varepsilon_{11} ,\; \varepsilon_{22} ,\; \gamma_{12} \right]^T$ and $\stressV := \left[\sigma_{11} ,\; \sigma_{22} ,\; \tau_{12} \right]^T$ and $\strainV := \left[\varepsilon_{11} ,\; \varepsilon_{22} ,\; \varepsilon_{33} ,\; \gamma_{12} ,\; \gamma_{13} \right]^T$ and $\stressV := \left[\sigma_{11} ,\; \sigma_{22} ,\; \sigma_{33} ,\; \tau_{12} ,\; \tau_{13} ,\; \tau_{23} \right]^T$ in two and three dimensions respectively, where the number of components of the strain and stress vectors is given by $\msd = \nsd(\nsd+1)/2$.

Using the Voigt notation, the relation between the displacement and the strain can be written as $\strainV = \gradS \bu$, where the matrix operator $\gradS \in \RR^{\msd \times \nsd}$ is given by
\begin{equation} \label{eq:symmGrad2D}
\gradS :=
\begin{bmatrix}
\partial/\partial x_1 & 0 & \partial/\partial x_2 \\
0 & \partial/\partial x_2 & \partial/\partial x_1
\end{bmatrix}^T
\end{equation}
and
\begin{equation} \label{eq:symmGrad3D}
\gradS :=
\begin{bmatrix}
\partial/\partial x_1 & 0 & 0 & \partial/\partial x_2 & \partial/\partial x_3 & 0 \\
0 & \partial/\partial x_2 & 0 & \partial/\partial x_1 & 0 & \partial/\partial x_3 \\
0 & 0 & \partial/\partial x_3 & 0 & \partial/\partial x_1 & \partial/\partial x_2
\end{bmatrix}^T,
\end{equation}
in two and three dimensions respectively.

The strain-stress relation given by Hooke's law also simplifies and can be written as $\stressV = \bD \strainV$, where $\bD \in \RR^{\msd \times \msd}$ is a symmetric and positive definite matrix that depends upon the material parameters characterising the medium and, in two dimensions, it also depends upon the model used (i.e. plane strain or plane stress). The matrix $\bD$ is given by
\begin{equation} \label{eq:LawVoigt2D}
\bD :=
\displaystyle\frac{E}{(1+\nu)(1-\vartheta\nu)}
\begin{bmatrix}
1+(1-\vartheta)\nu & \nu & 0 \\
\nu & 1+(1-\vartheta)\nu & 0 \\
0 & 0 & (1-\vartheta\nu)/2
\end{bmatrix}
\end{equation}
and
\begin{equation} \label{eq:LawVoigt3D}
\bD := \displaystyle\frac{E}{(1+\nu)(1-2\nu)}
\begin{bmatrix}
1-\nu & \nu & \nu & \\
\nu & 1-\nu & \nu & \bm{0}_{\nsd} \\
\nu & \nu & 1-\nu & \\
& \bm{0}_{\nsd} & & (1-2\nu)/2\Insd
\end{bmatrix},
\end{equation}
in two and three dimensions respectively, where $E$ is the Young modulus, $\nu$ the Poisson ratio. In two dimensions, the parameter $\vartheta=1$ denotes a plane stress model whereas $\vartheta=2$ denotes a plane strain model.

The strong form of the linear elastic problem can be written using Voigt notation as
\begin{equation} \label{eq:elasticitySystemVoigt}
\left\{\begin{aligned}
-\gradS^T \stressV  &= \bf       &&\text{in $\Omega$,}\\
\bu &= \bu_D  &&\text{on $\Gamma_D$,}\\
(1 - \xi) \Pn \bu +  (\Pt + \xi \Pn) \bN^T \stressV &= \xi \bg        &&\text{on $\Gamma_N$,}\\
\end{aligned}\right.
\end{equation}
where
\begin{equation} \label{eq:normalVoigt2D}
\bN :=
\begin{bmatrix}
n_1 & 0 & n_2 \\
0 & n_2 & n_1
\end{bmatrix}^T
\end{equation}
and
\begin{equation} \label{eq:normalVoigt3D}
\bN := \begin{bmatrix}
n_1 & 0 & 0 & n_2 & n_3 & 0\\
0 & n_2 & 0 & n_1 & 0 & n_3 \\
0 & 0 & n_3 & 0 & n_1 & n_2
\end{bmatrix}^T,
\end{equation}
in two and three dimensions respectively.

%==========================================================================
\section{Face centered finite volume (FCFV) formulation}
\label{sc:FCFV}
%==========================================================================

Let us introduce the \textit{broken computational domain} as a partition of the domain $\Omega$ in $\numel$ disjoint elements $\Omega_e$ with boundaries $\partial\Omega_e$.
The set of internal faces $\Gamma$ is defined as
\begin{equation} \label{eq:skeleton}
\Gamma := \left[ \bigcup_{e=1}^{\numel} \partial\Omega_e \right] \setminus \partial\Omega .
\end{equation}

In addition, the boundary of each element can be written as the union of its faces
\begin{equation}
\partial\Omega_e :=  \bigcup_{j=1}^{\numfa^e} \Gamma_{e,j},
\end{equation}
where $\numfa^e$ is the total number of faces of the element $\Omega_e$.

Following the standard notation using in HDG methods~\cite{RS-SH:16}, the discrete element spaces
\begin{subequations}
\begin{align}
\Vh(\Omega) & := \left\{ v \in \eltwo(\Omega) : v\vert_{\Omega_e} \in \Pk(\Omega_e) \forall \Omega_e \right\}, \\
\VhHat(S) & := \left\{ \hv \in \eltwo(S) : \hv\vert_{\Gamma_j}\in \Pk(\Gamma_j) \;\forall\Gamma_j\subset S\subseteq\Gamma\cup\Gamma_N \right\},
\end{align}
\end{subequations}
are introduced, where $\mathcal{P}^{k}(\Omega_e)$ and $\mathcal{P}^{k}(\Gamma_j)$ denote the space of polynomials of complete degree at most $k$ in $\Omega_e$ and on $\Gamma_j$ respectively.

%==========================================================================
\subsection{Mixed formulation}
\label{sc:mixedform}
%==========================================================================

In the proposed FCFV method, a mixed formulation of the elastic problem of Equation~\eqref{eq:elasticitySystemVoigt} is considered, namely
\begin{equation} \label{eq:elasticityBrokenFirstOrder}
\left\{\begin{aligned}
\bL + \bDHalf \gradS \bu &= \bm{0}    &&\text{in $\Omega_e$, and for $e=1,\dotsc ,\numel$,}\\	
\gradS^T \bDHalf \bL &= \bm{f}          &&\text{in $\Omega_e$, and for $e=1,\dotsc ,\numel$,}\\
\bu &= \bu_D     &&\text{on $\Gamma_D$,}\\
(1 - \xi) \Pn \bu -  (\Pt + \xi \Pn) \bN^T \bDHalf \bL  &= \xi \bg         &&\text{on $\Gamma_N$,}\\
\jump{\bu \otimes \bn} &= \bm{0}  &&\text{on $\Gamma$,}\\
\jump{\bN^T \bDHalf \bL} &= \bm{0}  &&\text{on $\Gamma$,}\\
\end{aligned} \right.
\end{equation}
where the \emph{jump} operator $\jump{\cdot}$ is defined over an internal face shared by two elements $\Omega_e$ and $\Omega_l$ as the sum of the values from the two elements sharing the face~\cite{AdM-MFH:08}, namely
\begin{equation}
\jump{\odot} = \odot_e + \odot_l.
\end{equation}

The matrix $\bDHalf$, introduced in Equation~\eqref{eq:elasticityBrokenFirstOrder} to guarantee the symmetry of the mixed formulation, is defined as $\bDHalf = \bV \bLambda^{1/2} \bV^T$, after performing the spectral decomposition of the matrix $\bD = \bV \bLambda \bV^T$, where the matrix $\bV$ and the diagonal matrix $\bLambda$ contain the eigenvectors and eigenvalues of $\bD$ respectively and $\bLambda^{1/2}$ is the diagonal matrix containing the square root of the eigenvalues of $\bD$.

It is worth noting that the last two equations in~\eqref{eq:elasticityBrokenFirstOrder}, called \textit{transmission conditions}, impose the continuity of the displacement field and the normal stress across the internal faces $\Gamma$.

%==========================================================================
\subsection{FCFV weak formulation}
\label{sc:FCFVweakform}
%==========================================================================

As other HDG methods~\cite{Jay-CGL:09,Cockburn-CDG:08,Nguyen-NPC:09,Nguyen-NPC:09b,Nguyen-NPC:10,Nguyen-NPC:11,RS-SH:16,RS-AH:18}, the proposed FCFV method solves the mixed problem in the broken computational domain in two phases. First, a purely Dirichlet problem is defined on each element to write the mixed and primal variables $(\bL_e,\bu_e)$ in terms of a new hybrid variable $\bhu$, uniquely defined as the trace of the displacement field on $\Gamma \cup \Gamma_N$, namely
\begin{equation} \label{eq:elasticityStrongLocal}
\left\{\begin{aligned}
\bL_e + \bDHalf \gradS \bu_e &= \bm{0}    &&\text{in $\Omega_e$}\\	
\gradS^T \bDHalf \bL_e &= \bm{f}          &&\text{in $\Omega_e$}\\
\bu_e &= \bu_D     &&\text{on $\partial\Omega_e \cap \Gamma_D$,}\\
\bu_e &= \bhu  &&\text{on $\partial\Omega_e \setminus \Gamma_D$,}\\
\end{aligned} \right.
\end{equation}
for $e=1,\dotsc ,\numel$. This set of problems are usually referred to as \textit{local problems} and the solution in one element is independent on the solution in the other elements.

Second, the so-called \textit{global problem} is defined to compute the hybrid variable $\bhu$, namely
\begin{equation} \label{eq:elasticityStrongGlobalPrevious}
\left\{\begin{aligned}
\jump{\bu \otimes \bn} &= \bm{0}  &&\text{on $\Gamma$,}\\
\jump{\bN^T \bDHalf \bL} &= \bm{0}  &&\text{on $\Gamma$,}\\
(1 - \xi) \Pn \bu -  (\Pt + \xi \Pn) \bN^T \bDHalf \bL  &= \xi \bg  &&\text{on $\Gamma_N$.}\\
\end{aligned} \right.
\end{equation}
It is worth noting that the first transmission condition in Equation~\eqref{eq:elasticityStrongGlobalPrevious} is automatically satisfied due to the imposed Dirichlet boundary conditions in the local problems of Equation~\eqref{eq:elasticityStrongLocal} and the unique definition of the hybrid variable on each face. Therefore, the global problem is simply
\begin{equation} \label{eq:elasticityStrongGlobal}
\left\{\begin{aligned}
\jump{\bN^T \bDHalf \bL} &= \bm{0}  &&\text{on $\Gamma$,}\\
(1 - \xi) \Pn \bu -  (\Pt + \xi \Pn) \bN^T \bDHalf \bL  &= \xi \bg  &&\text{on $\Gamma_N$.}\\
\end{aligned} \right.
\end{equation}

Next, the weak formulation of both the local and global problems is presented. For each element $\Omega_e, \ e=1,\dotsc ,\numel$, the weak formulation of Equation~\eqref{eq:elasticityStrongLocal} reads as follows: given $\bu_D$ on $\Gamma_D$ and $\bhu $ on $\Gamma\cup\Gamma_N$, find $(\bL^h_e ,\bu^h_e) \in [\Vh(\Omega_e)]^{\msd} \times [\Vh(\Omega_e)]^{\nsd}$ such that
\begin{subequations}\label{eq:HDGElasticityWeakLocalPre}
	\begin{align}
	&
	- (\bv,\bL^h_e)_{\Omega_e} + (\gradS^T \bDHalf \bv, \bu^h_e)_{\Omega_e} =   \langle \bN_e^T \bDHalf \bv , \bu_D\rangle_{\partial\Omega_e\cap\Gamma_D} + \langle \bN_e^T \bDHalf \bv , \bhu^h \rangle_{\partial\Omega_e\setminus\Gamma_D} , \label{eq:HDGElasticityWeakLocalLPre}
	\\
	&
	-(\gradS \bw, \bDHalf \bL^h_e)_{\Omega_e}  + \langle \bw ,\bN_e^T \widehat{\bDHalf \bL^h_e} \rangle_{\partial\Omega_e} =  (\bw,\bm{f})_{\Omega_e} , \label{eq:HDGElasticityWeakLocalUPre}
	\end{align}
\end{subequations}
for all $(\bv ,\bw) \in [\Vh(\Omega_e)]^{\msd} \times [\Vh(\Omega_e)]^{\nsd}$.

In the above expressions, the following definition of the internal products of vector functions in $\eltwo(\Omega_e)$ has been used:
\begin{equation} \label{eq:innerScalar}
(\bm{p},\bm{q})_{\Omega_e} := \int_{\Omega_e} \bm{p} \cdot\bm{q} \ d\Omega  , \qquad \langle \hat{\bm{p}}, \hat{\bm{q}} \rangle_{\partial\Omega_e} := \sum_{j=1}^{\numfa^e}  \int_{\Gamma_{e,j}} \hat{\bm{p}} \cdot \hat{\bm{q}} \ d\Gamma.
\end{equation}

As usual in an HDG context, Dirichlet boundary conditions are imposed in the weak form and the trace of the numerical stress is defined as
\begin{equation} \label{eq:traceElasticity}
\bN_e^T \widehat{\bDHalf \bL^h_e} := 
\begin{cases}
\bN_e^T \bDHalf \bL^h_e + \btau_e (\bu^h_e - \bu_D) & \text{on $\partial\Omega_e\cap\Gamma_D$,} \\
\bN_e^T \bDHalf \bL^h_e + \btau_e (\bu^h_e - \bhu^h) & \text{elsewhere,}  
\end{cases}
\end{equation}
where the stabilisation tensor $\btau_e$ is introduced to ensure the stability, accuracy and convergence of the resulting numerical scheme~\cite{Jay-CGL:09,Cockburn-CDG:08}.

Integrating again by parts Equation~\eqref{eq:HDGElasticityWeakLocalUPre} and introducing the definition of Equation~\eqref{eq:traceElasticity}, the weak formulation of the local problems, for $e=1,\dotsc ,\numel$, reads: given $\bu_D$ on $\Gamma_D$ and $\bhu $ on $\Gamma\cup\Gamma_N$, find $(\bL^h_e ,\bu^h_e) \in [\Vh(\Omega_e)]^{\msd} \times [\Vh(\Omega_e)]^{\nsd}$ such that
\begin{subequations}\label{eq:HDGElasticityWeakLocal}
\begin{align} 
&	- (\bv,\bL^h_e)_{\Omega_e} + (\gradS^T \bDHalf \bv, \bu^h_e)_{\Omega_e} =   \langle \bN_e^T \bDHalf \bv , \bu_D\rangle_{\partial\Omega_e\cap\Gamma_D} + \langle \bN_e^T \bDHalf \bv , \bhu^h \rangle_{\partial\Omega_e\setminus\Gamma_D} , \label{eq:HDGElasticityWeakLocalL}
	\\	
&	(\bw, \gradS^T \bDHalf \bL^h_e)_{\Omega_e}  + \langle \bw , \btau_e \bu^h_e \rangle_{\partial\Omega_e} =  (\bw,\bm{f})_{\Omega_e}  + \langle \bw , \btau_e \bu_D\rangle_{\partial\Omega_e\cap\Gamma_D} + \langle \bw , \btau_e \bhu^h \rangle_{\partial\Omega_e\setminus\Gamma_D} , \label{eq:HDGElasticityWeakLocalU}
\end{align}
\end{subequations}
for all $(\bv ,\bw) \in [\Vh(\Omega_e)]^{\msd} \times [\Vh(\Omega_e)]^{\nsd}$.

Analogously, the weak formulation of the global problem is found by multiplying by a vector of test functions in $[\VhHat(\Gamma\cup\Gamma_N)]^{\nsd}$ and adding all the contributions corresponding to internal faces and faces on the Neumann boundary. It reads, find $\bhu^h\in[\VhHat(\Gamma\cup\Gamma_N)]^{\nsd}$ that satisfies
\begin{multline} \label{eq:HDGElasticityWeakGlobalPrevious}
\sum_{e=1}^{\numel}\left\{
\langle \bhw, \bN_e^T \widehat{\bDHalf \bL^h_e} \rangle_{\partial\Omega_e\setminus\partial\Omega}
+ \langle \bhw, (1 - \xi) \Pn \bhu^h -  (\Pt + \xi \Pn) \bN^T \widehat{\bDHalf \bL^h_e} \rangle_{\partial\Omega_e\cap\Gamma_N} \right\} \\
 = \sum_{e=1}^{\numel} \langle \bhw, \xi \bg \rangle_{\partial\Omega_e\cap\Gamma_N},
\end{multline}
for all $\bhw\in[\VhHat(\Gamma\cup\Gamma_N)]^{\nsd}$.

By introducing the definition of Equation~\eqref{eq:traceElasticity}, the weak formulation of the global problem is: find $\bhu^h\in[\VhHat(\Gamma\cup\Gamma_N)]^{\nsd}$ such that
\begin{multline} \label{eq:HDGElasticityWeakGlobal}
\sum_{e=1}^{\numel}\left\{
\langle \bhw, \bN_e^T \bDHalf \bL^h_e \rangle_{\partial\Omega_e\setminus\partial\Omega}
+ \langle \bhw,\btau_e\, \bu^h_e \rangle_{\partial\Omega_e\setminus\partial\Omega} 
- \langle \bhw,\btau_e\,\bhu^h \rangle_{\partial\Omega_e\setminus\partial\Omega}\right. \\
-\langle \bhw, (\Pt + \xi \Pn) \bN_e^T \bDHalf \bL^h_e \rangle_{\partial\Omega_e\cap\Gamma_N}
- \langle \bhw, (\Pt + \xi \Pn) \btau_e\, \bu^h_e \rangle_{\partial\Omega_e\cap\Gamma_N} 
\\
\left.
+ \langle \bhw, \left[ (1 - \xi) \Pn + (\Pt + \xi \Pn) \btau_e \right]\,\bhu^h \rangle_{\partial\Omega_e\cap\Gamma_N} \right\}
= \sum_{e=1}^{\numel} \langle \bhw, \xi \bg \rangle_{\partial\Omega_e\cap\Gamma_N},
\end{multline}
for all $\bhw\in[\VhHat(\Gamma\cup\Gamma_N)]^{\nsd}$.

%==========================================================================
\subsection{FCFV spatial discretisation}
\label{sc:HDGdiscretisation}
%==========================================================================

The proposed FCFV rationale consists of employing a constant degree of approximation within each element for the mixed and primal variables $\bL_e$ and $\bu_e$ and a constant degree of approximation on each face for the hybrid variable $\bhu$. The discrete form of the local problem of Equation~\eqref{eq:HDGElasticityWeakLocal} is obtained as
\begin{subequations}\label{eq:HDGElasticityWeakLocalK0}
\begin{align}
&
-| \Omega_e | \Le =  \sum_{j \in \Dset} | \Gamma_{e,j} | \bDHalf_j^T \bN_j \bu_{D,j}   + \sum_{j \in \Bset} | \Gamma_{e,j} | \bDHalf_j^T \bN_j  \buHj  , \label{eq:HDGElasticityWeakLocalLK0}
\\
&
\sum_{j \in \Aset} | \Gamma_{e,j} | \btau_j \bue =  | \Omega_e | \bm{f}_e  + \sum_{j \in \Dset} | \Gamma_{e,j} | \btau_j \bu_{D,j} + \sum_{j \in \Bset} | \Gamma_{e,j} | \btau_j \buHj, \label{eq:HDGElasticityWeakLocalUK0}
\end{align}
\end{subequations}
for $e=1,\dotsc,\numel$, where $\Le$ and $\bue$ denote the constant value of the mixed and primal variables in the element $\Omega_e$, $\buHj$ denotes the constant value of the hybrid variable on the face $\Gamma_{e,j}$ and the following sets of faces have been introduced for each element:
\begin{equation}\label{eq:setsFaces}
\Aset :=  \{1, \ldots, \numfa^e \},  \qquad
\Dset := \{j \in \Aset \; | \; \Gamma_{e,j} \cap \Gamma_D \neq \emptyset \}, \qquad
\Bset := \Aset \setminus \Dset ,
\end{equation}
with $\numfa^e$ the total number of faces of $\Omega_e$.

It is worth noting that the discrete form of the local problem has been obtained by utilising a quadrature with one integration point to compute the integrals of the weak formulation.

An important advantage of using of a constant degree of approximation for the mixed and primal variables is that the two equations of the local problem decouple and it is possible to obtain a closed form expression for $\Le$ and $\bue$ as a function of $\buHj$, namely
\begin{subequations}\label{eq:HDGElasticityWeakLocalK0b}
\begin{align}
&
\Le =  -| \Omega_e |^{-1} \vect{z}_e   - | \Omega_e |^{-1} \sum_{j \in \Bset} | \Gamma_{e,j} | \bDHalf_j^T \bN_j \buHj  , \label{eq:HDGElasticityWeakLocalLK0b}
\\
&
\bue =  \bm{\alpha}_e^{-1} \bm{\beta}_e + \bm{\alpha}_e^{-1} \sum_{j \in \Bset} | \Gamma_{e,j} | \btau_j \buHj, \label{eq:HDGElasticityWeakLocalUK0b}
\end{align}
\end{subequations}
where 
\begin{equation} \label{eq:stokesPrecomp}
\begin{split}
\bm{\alpha}_e  := \sum_{j \in \Aset} | \Gamma_{e,j} & | \btau_j , 
\qquad
\bm{\beta}_e  :=  | \Omega_e | \bm{f}_e  + \sum_{j \in \Dset} | \Gamma_{e,j} | \btau_j \bu_{D,j},  
\\
& \vect{z}_e  := \sum_{j \in \Dset} | \Gamma_{e,j} | \bDHalf_j^T \bN_j \bu_{D,j}.
\end{split}
\end{equation}

Similarly, employing a constant degree of approximation for $\bL_e$, $\bu_e$ and $\bhu$, the discrete form of the global problem of Equation~\eqref{eq:HDGElasticityWeakGlobal} is 
\begin{equation} \label{eq:HDGElasticityWeakGlobalK0}
\sum_{e=1}^{\numel}
| \Gamma_{e,i} |
\left\{
\mat{A}_{e,i} \bN_i^T \bDHalf_i \Le + \mat{A}_{e,i}  \btau_i \bue +  \mat{B}_{e,i} \buHi
\right\} 
= \xi \sum_{e=1}^{\numel} | \Gamma_{e,i} | \bg_i \, \chi_{\Nset}(i),
\end{equation}
for $i \in \Bset$, where $\chi_{\Iset}$ and $\chi_{\Nset}$ are the indicator functions of the sets $\Iset := \{j \in \Aset \; | \; \Gamma_{e,j} \cap \partial \Omega = \emptyset \}$ and $\Nset := \{j \in \Aset \; | \; \Gamma_{e,j} \cap \Gamma_N \neq \emptyset \}$ respectively. The following matrices have been introduced to shorten the notation in Equation~\eqref{eq:HDGElasticityWeakGlobalK0}
\begin{subequations}
	\begin{align}
	&
	\mat{A}_{e,i} =  \Insd \chi_{\Iset}(i) -  (\Pt + \xi \Pn) \chi_{\Nset}(i)  ,
	\\
	&
	\mat{B}_{e,i} =  -\btau_i \chi_{\Iset}(i) +  \left[ (1 - \xi) \Pn + (\Pt + \xi \Pn) \btau_i \right] \chi_{\Nset}(i) .
	\end{align}
\end{subequations}

After introducing the closed form expressions of the mixed and primal variable of Equation~\eqref{eq:HDGElasticityWeakLocalK0b} in Equation~\eqref{eq:HDGElasticityWeakGlobalK0}, a linear system of equations is obtained, where the only unknown is the hybrid variable defined over the interior and Neumann faces, $\Gamma \cup \Gamma_N$, namely
\begin{equation} \label{eq:globalSystemElasticityK0}
\mat{\widehat{K}} \vect{\hat{u}}  = \vect{\hat{f}}.
\end{equation}

The matrix $\mat{\widehat{K}}$ and the vector $\vect{\hat{f}}$ are obtained by assembling the elemental contributions
\begin{subequations}\label{HDG-Elasticity-globalSystemK0}
\begin{align}
\mat{\widehat{K}}^e_{i,j} & :=   | \Gamma_{e,i} | \left\{  |\Gamma_{e,j}| \mat{A}_{e,i} \bigl( \btau_i \bm{\alpha}_e^{-1} \btau_j -  |\Omega_e|^{-1} \bN_i^T \bDHalf_i \bDHalf_j^T \bN_j \bigr) + \mat{B}_{e,i} \delta_{ij} \right\}
, \\
\vect{\widehat{f}}^e_i & :=   | \Gamma_{e,i} | \left\{ - \mat{A}_{e,i} \bigl( \btau_i \bm{\alpha}_e^{-1} \bm{\beta}_e - |\Omega_e|^{-1} \bN_i^T \bDHalf_i \vect{z}_e \bigr) + \xi \bg_i \, \chi_{\Nset}(i) \right\} ,
\end{align}
\end{subequations}
for $i,j \in \Bset$ and $\delta_{ij}$ being the classical Kronecker delta, equal to $1$ if $i=j$ and $0$ otherwise.

%==========================================================================
\section{Two dimensional examples}
\label{sc:examples2D}
%==========================================================================

This Section presents three numerical examples in two dimensions. The first example is used in order to validate the optimal rate of convergence, to illustrate the robustness of the proposed FCFV approach in the incompressible limit, to numerically study the effect of the stabilisation parameter and to demonstrate the robustness in terms of element distortion. The last two examples involve classical test cases for linear elastic solvers, namely the Kirsch's plate and the Cook's membrane problems.

%==========================================================================
\subsection{Optimal order of convergence}
\label{sc:convergence2D}
%==========================================================================

The first example considers a mesh convergence study to verify the optimal approximation properties of the proposed FCFV method in two dimensions. The model problem of Equation~\eqref{eq:elasticity}, defined in $\Omega = [0,1]^2$, is considered. The external force and boundary conditions are selected so that the exact solution~\cite{soon2009hybridizable} is given by 
\begin{subequations}
\begin{align}
u_1(x_1,x_2) & = -x_1^2 x_2 (x_1-1)^2(x_2-1)(2x_2-1), \\
u_2(x_1,x_2) & =  x_2^2 x_1 (x_2-1)^2(x_1-1)(2x_1-1).
\end{align}
\end{subequations}
The traction corresponding to the analytical solution is imposed on $\Gamma_N = \{(x_1,x_2) \in \mathbb{R}^2 \; | \; x_2=0\}$, whereas homogeneous Dirichlet boundary conditions are imposed on the rest of the boundary.

Structured uniform quadrilateral meshes with characteristic element size $h=2^{-r}$ are generated, where $r$ denotes the level of mesh refinement. Triangular uniform meshes are obtained by subdivision of each quadrilateral in four triangles using the two diagonals of the quadrilateral.

\begin{figure}[!tb]
	\centering
	\subfigure[Mesh 3]{\includegraphics[width=0.32\textwidth]{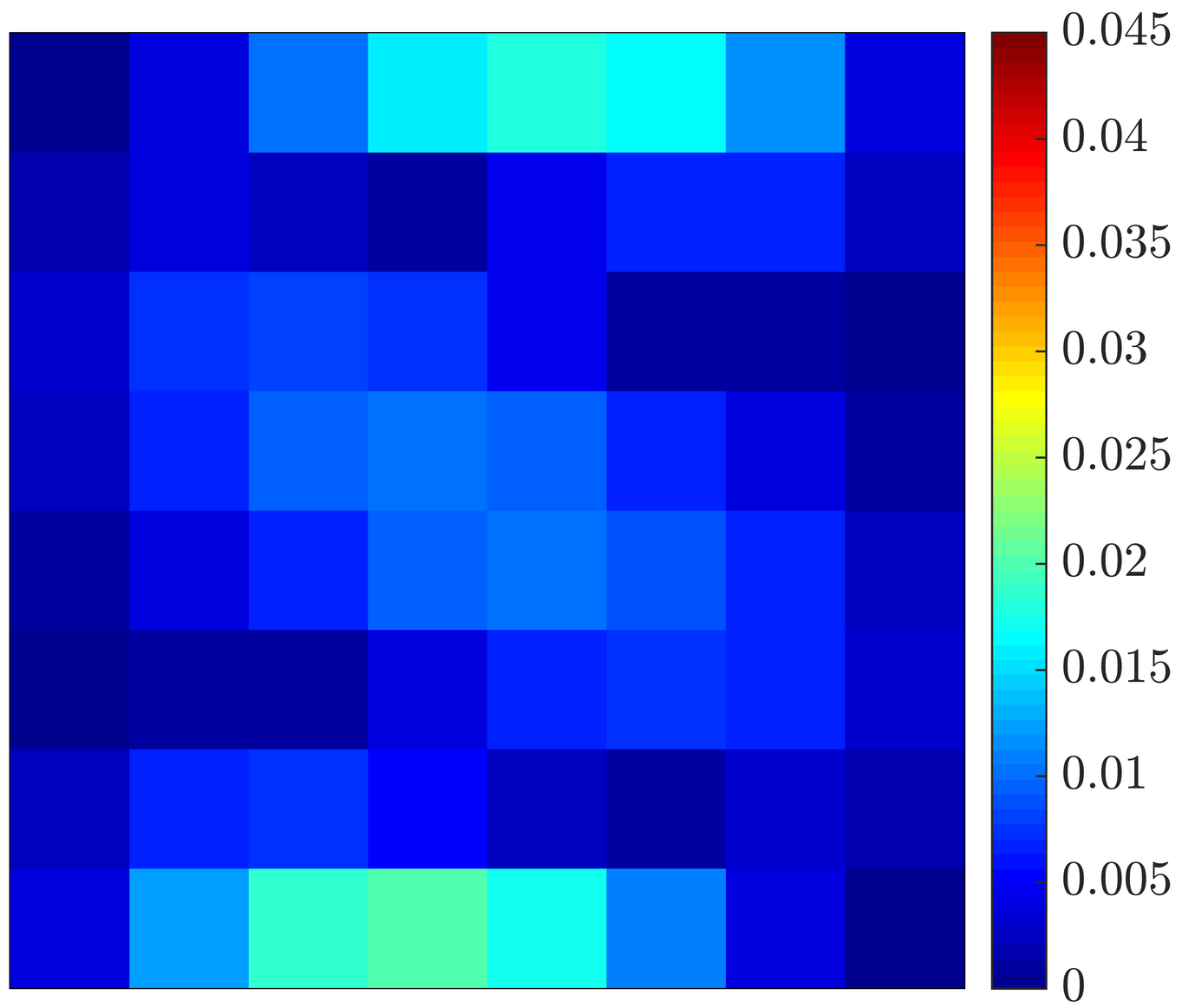}}
	\subfigure[Mesh 5]{\includegraphics[width=0.32\textwidth]{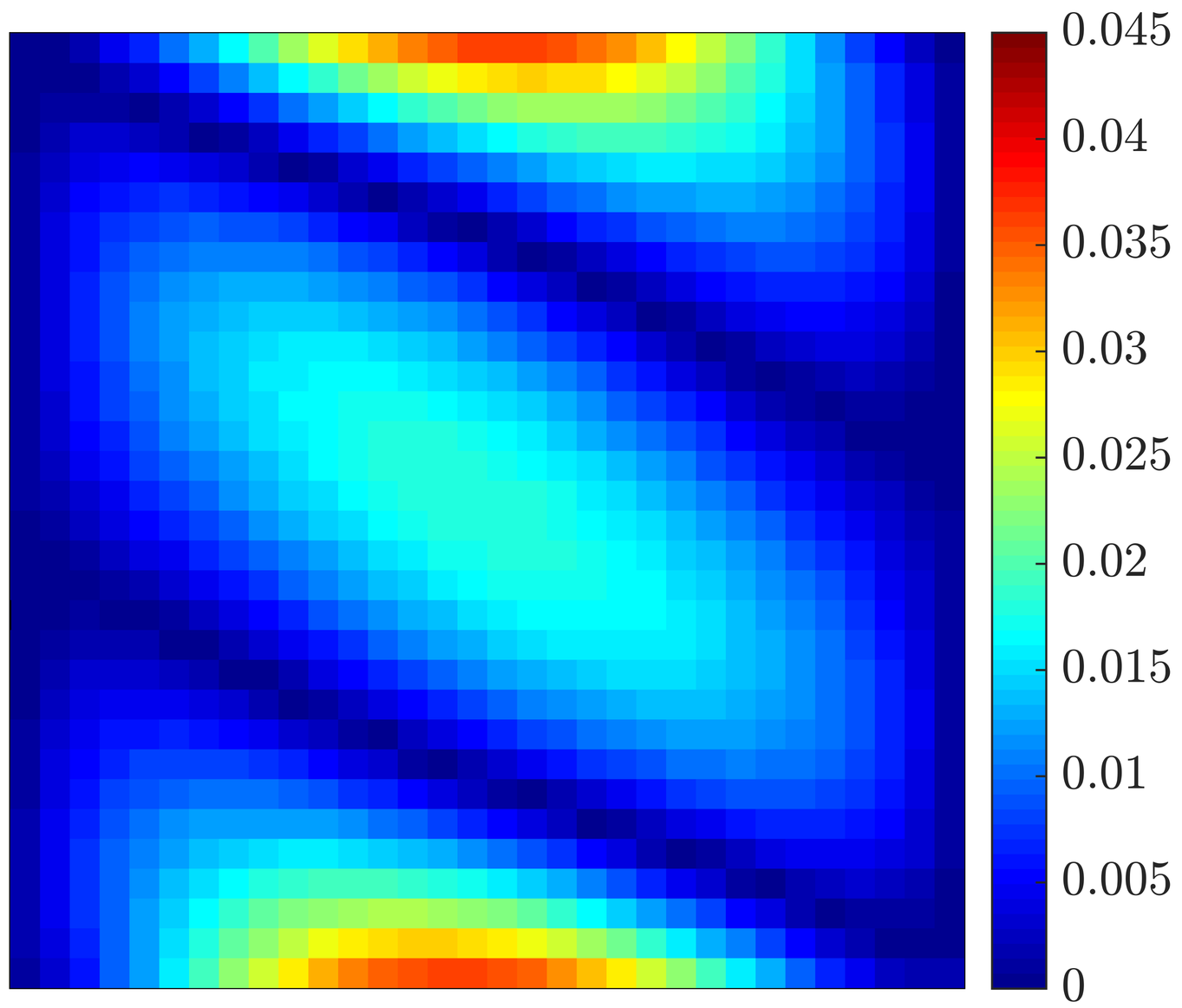}}
	\subfigure[Mesh 7]{\includegraphics[width=0.32\textwidth]{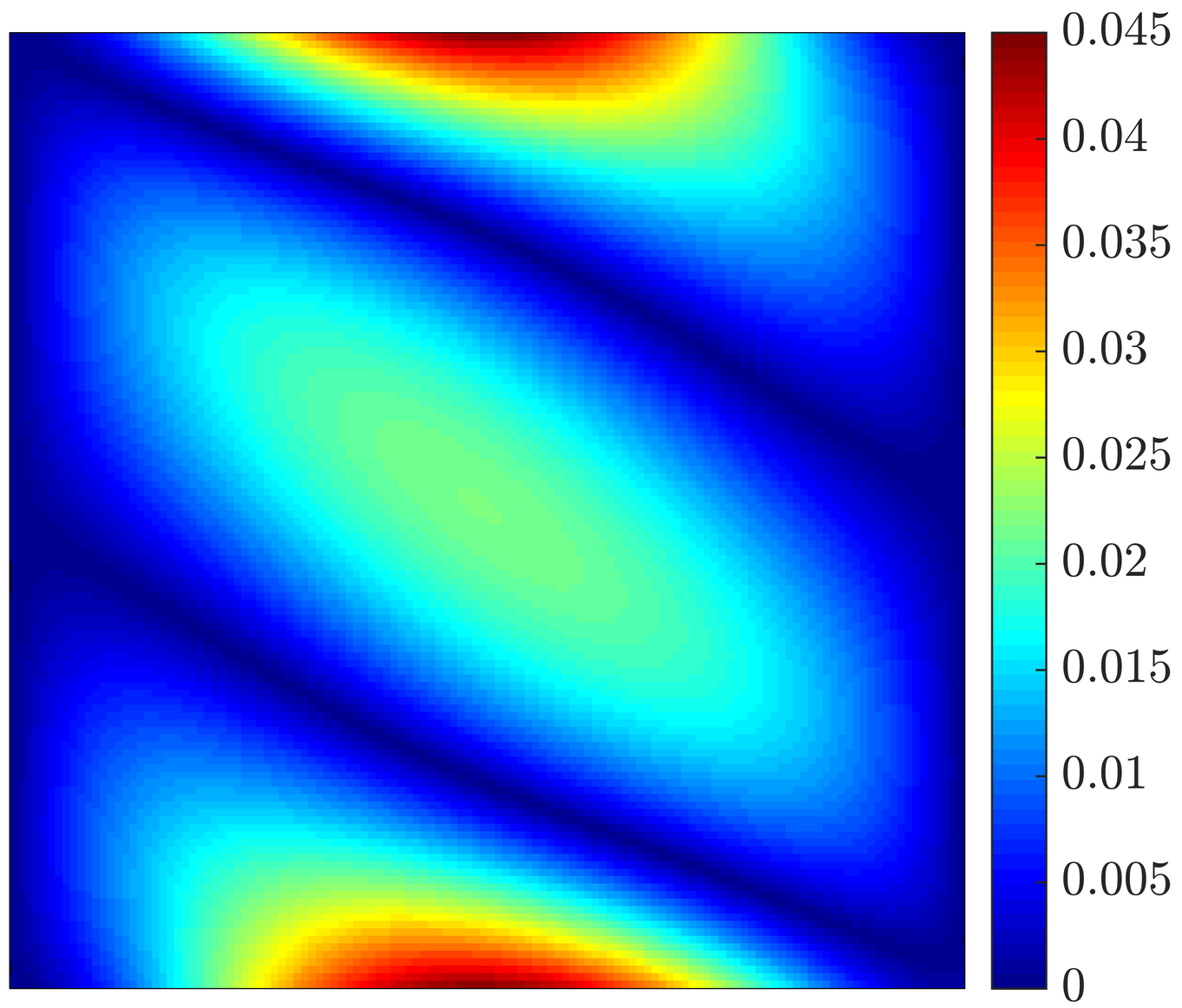}}
	\caption{Von Misses stress on three quadrilateral meshes.}
	\label{fig:2DsolVM}
\end{figure}

The computed Von Misses stress on three quadrilateral meshes is represented in Figure~\ref{fig:2DsolVM}, illustrating the increasing accuracy offered by the proposed FCFV as the mesh is refined.

\begin{figure}[!tb]
	\centering
	\subfigure[Displacement]{\includegraphics[width=0.49\textwidth]{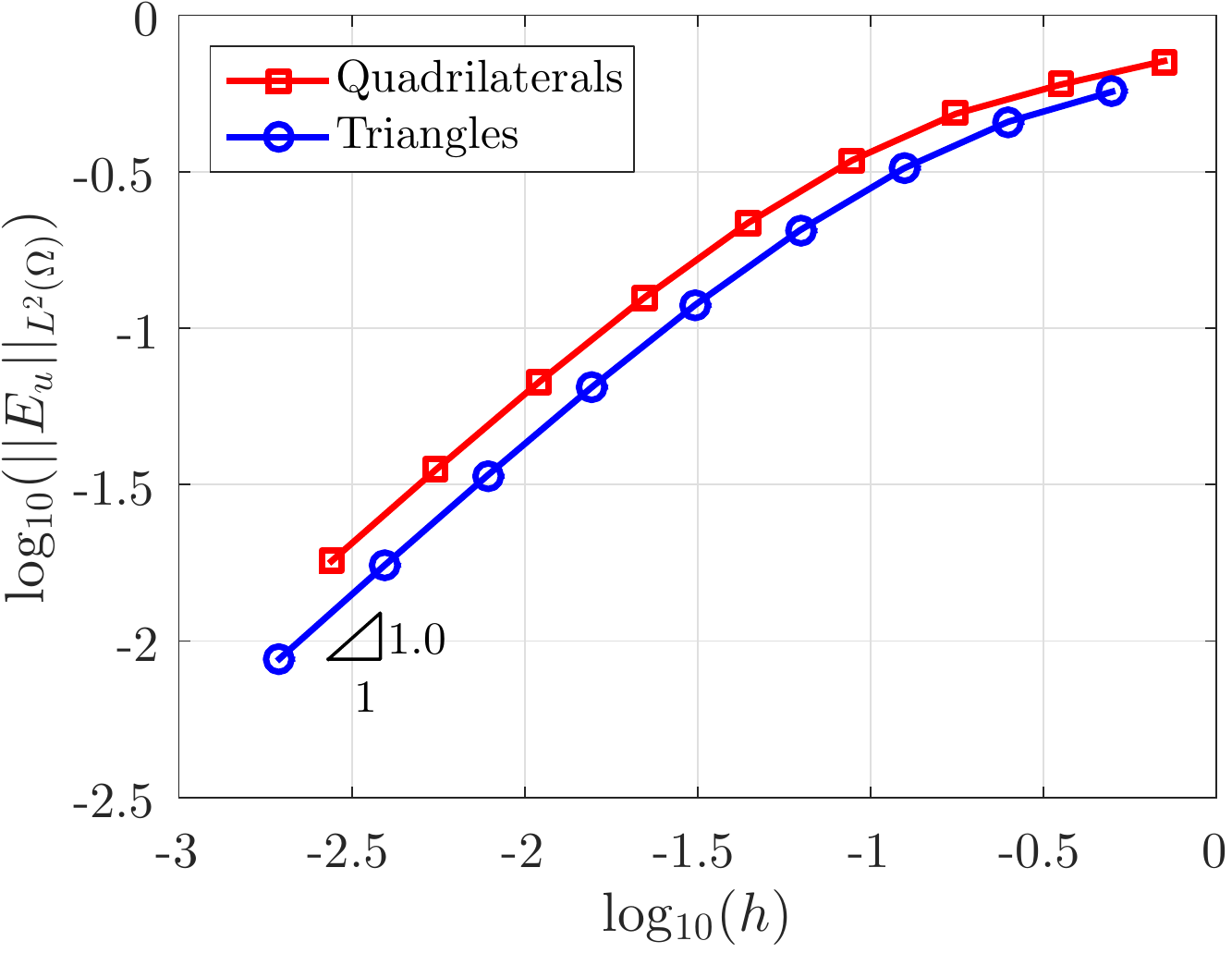}}
	\subfigure[Stress]{\includegraphics[width=0.49\textwidth]{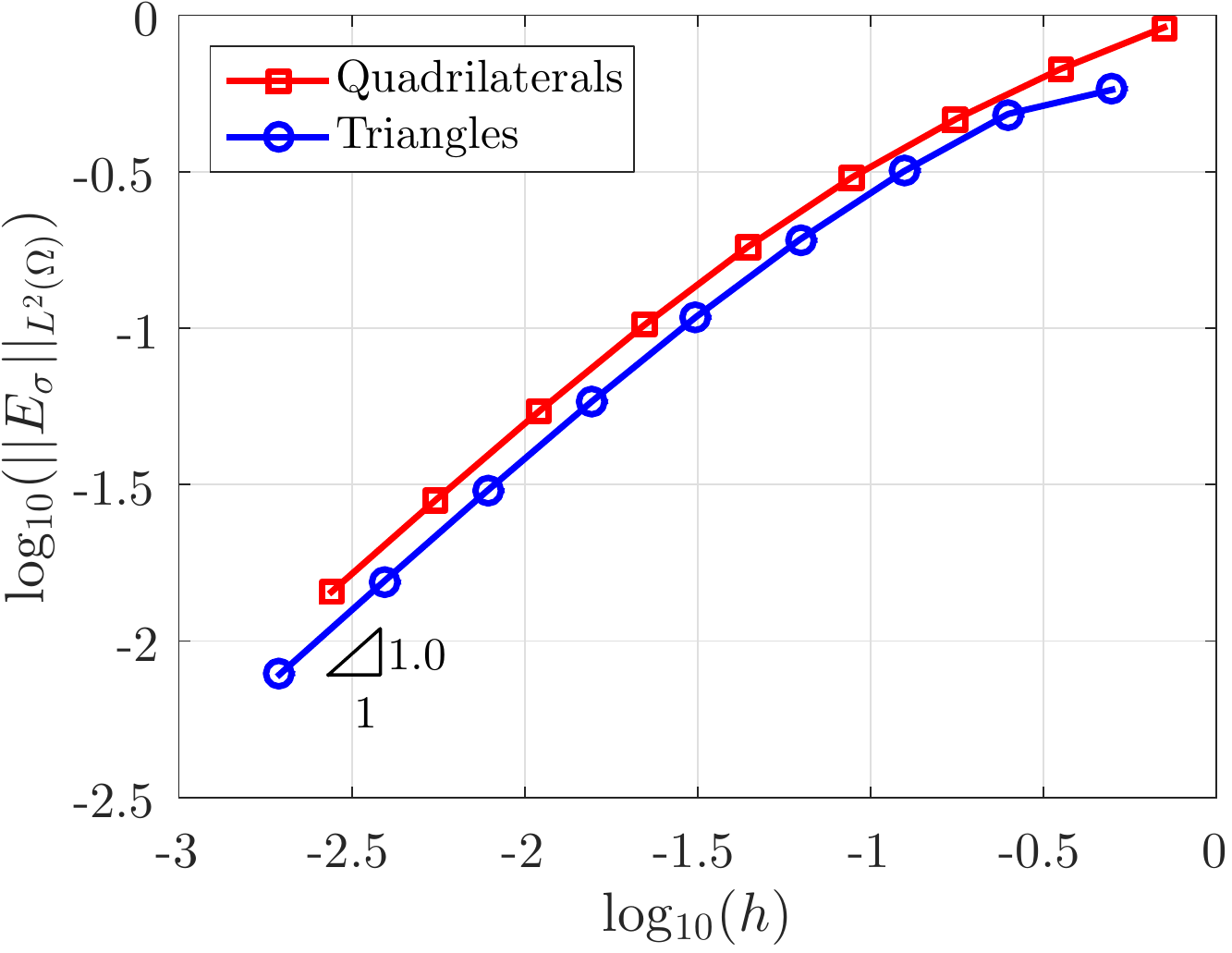}}	
	\caption{Mesh convergence of the $\eltwo(\Omega)$ error of the (a) displacement and (b) the stress, for quadrilateral and triangular elements in a medium with $\nu=1/3$.}
	\label{fig:convergence2D_nu0p3333}
\end{figure}

Figure~\ref{fig:convergence2D_nu0p3333} displays the error of the computed displacement and the stress fields in the $\eltwo(\Omega)$ norm as a function of the characteristic element size for both quadrilateral and triangular meshes on a medium with $E=1$ and $\nu=1/3$. 
It can be observed that the error converges with the expected rate of convergence for both the displacement and the stress. It is important to emphasise that the proposed FCFV produces a stress field with an error that converges linearly to the exact solution without performing a reconstruction of the displacement field. In addition, the proposed approach provides similar accuracy for both the displacement and the stress field due to the use of a mixed formulation.

%==========================================================================
\subsection{Locking-free behaviour for nearly incompressible materials}
\label{sc:incompressible}
%==========================================================================

\begin{figure}[!tb]
	\centering
	\subfigure[Displacement]{\includegraphics[width=0.49\textwidth]{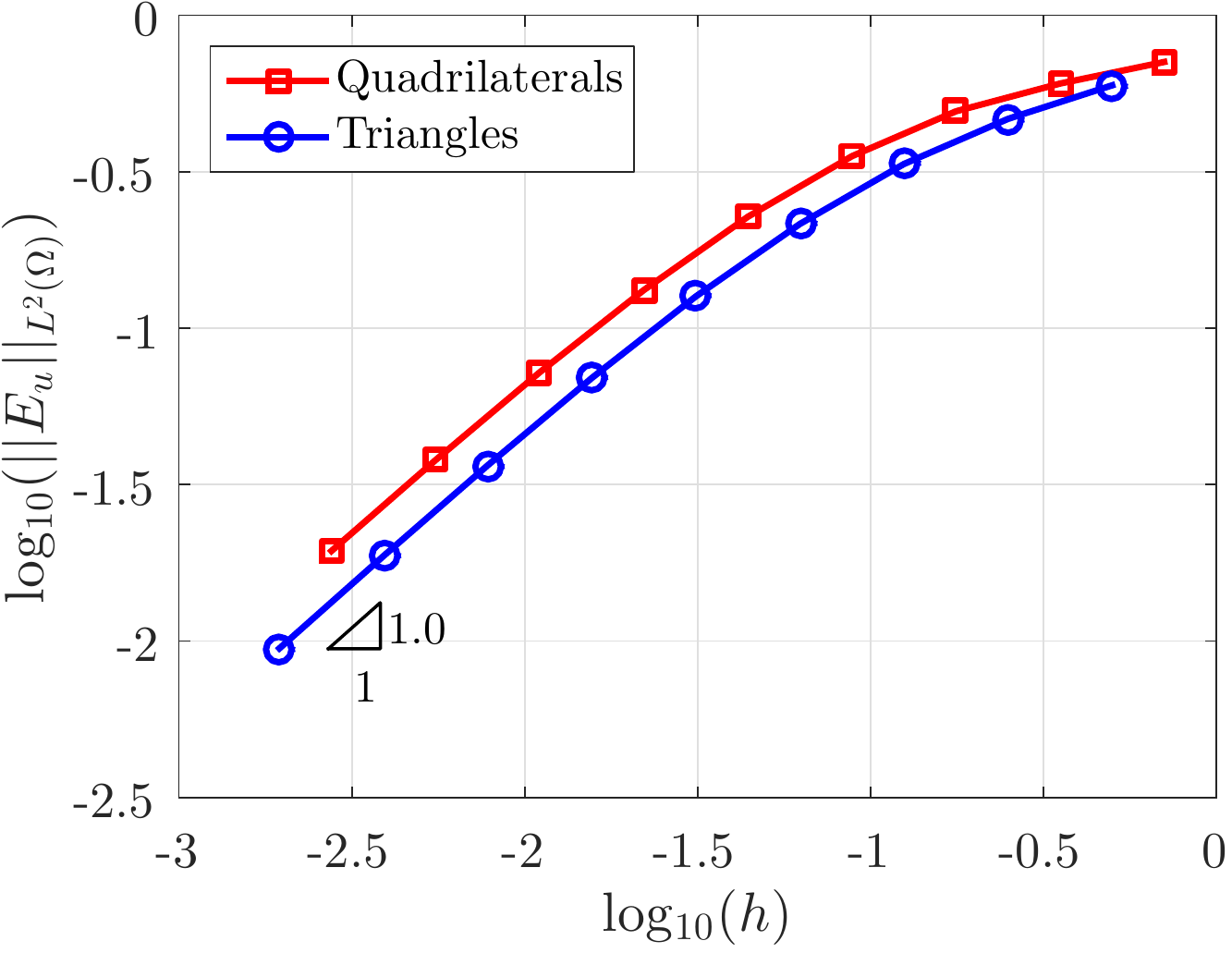}}
	\subfigure[Stress]{\includegraphics[width=0.49\textwidth]{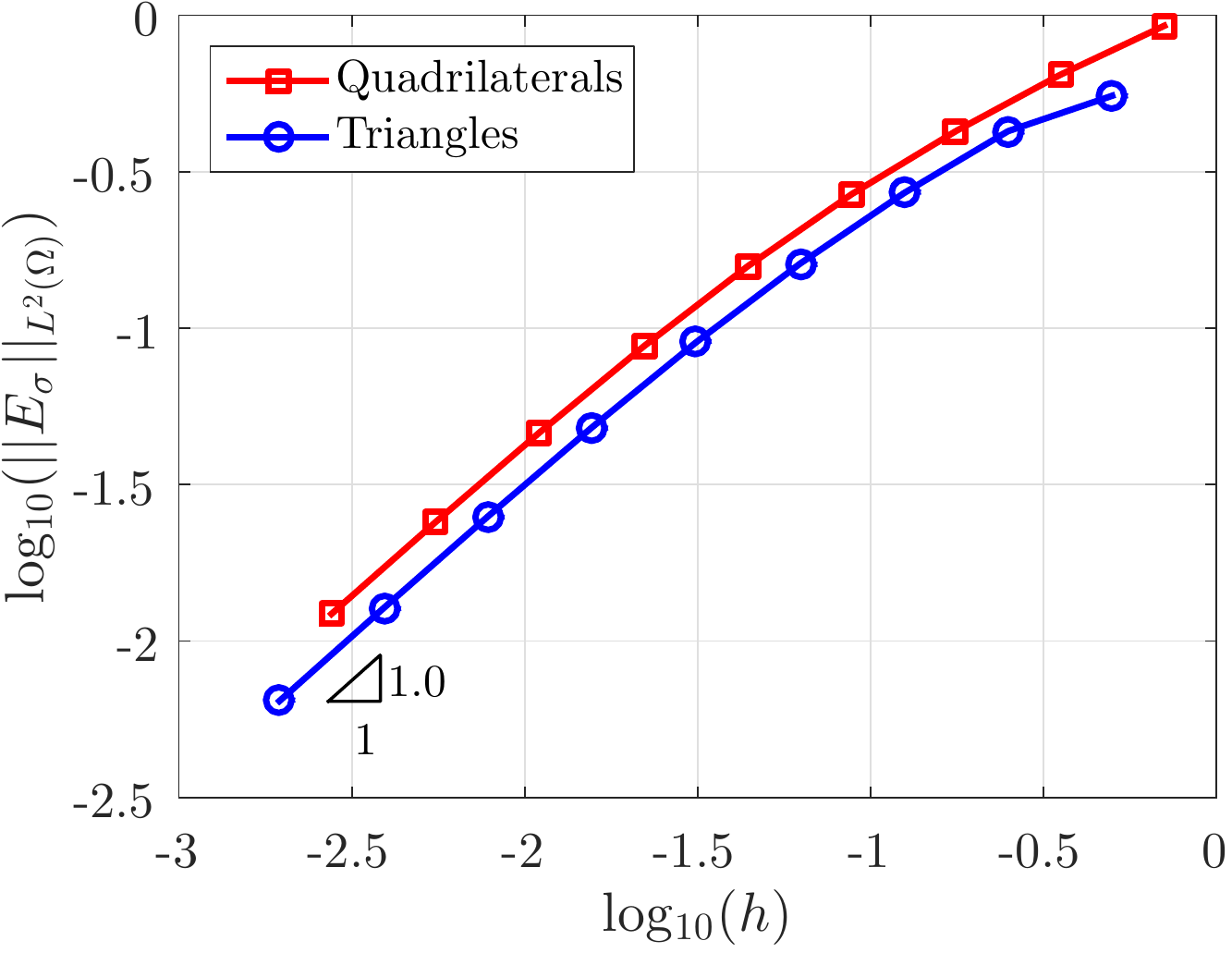}}	
	\caption{Mesh convergence of the $\eltwo(\Omega)$ error of the (a) displacement and (b) the stress, for quadrilateral and triangular elements in a medium with $\nu=0.49999$.}
	\label{fig:convergence2D_nu0p49999}
\end{figure}

To demonstrate the robustness of the proposed approach for nearly incompressible materials, the problem considered in Section~\ref{sc:convergence2D} is studied for a material with $E=1$ and $\nu=0.49999$. Figure~\ref{fig:convergence2D_nu0p49999} displays the error of the computed displacement and the stress fields in the $\eltwo(\Omega)$ norm as a function of the characteristic element size for both quadrilateral and triangular meshes. 
The results exhibit the optimal order of convergence for both the displacement and the stress. In addition, by comparing Figures~\ref{fig:convergence2D_nu0p49999} and \ref{fig:convergence2D_nu0p3333} it can be observed that almost identical results are obtained irrespectively of the value of the Poisson ratio, illustrating the robustness and suitability of the proposed approach for nearly incompressible materials.

%==========================================================================
\subsection{Optimal value of the stabilisation parameter}
\label{sc:tau2D}
%==========================================================================

\begin{figure}[!tb]
	\centering
	\subfigure[Quadrilaterals]{\includegraphics[width=0.49\textwidth]{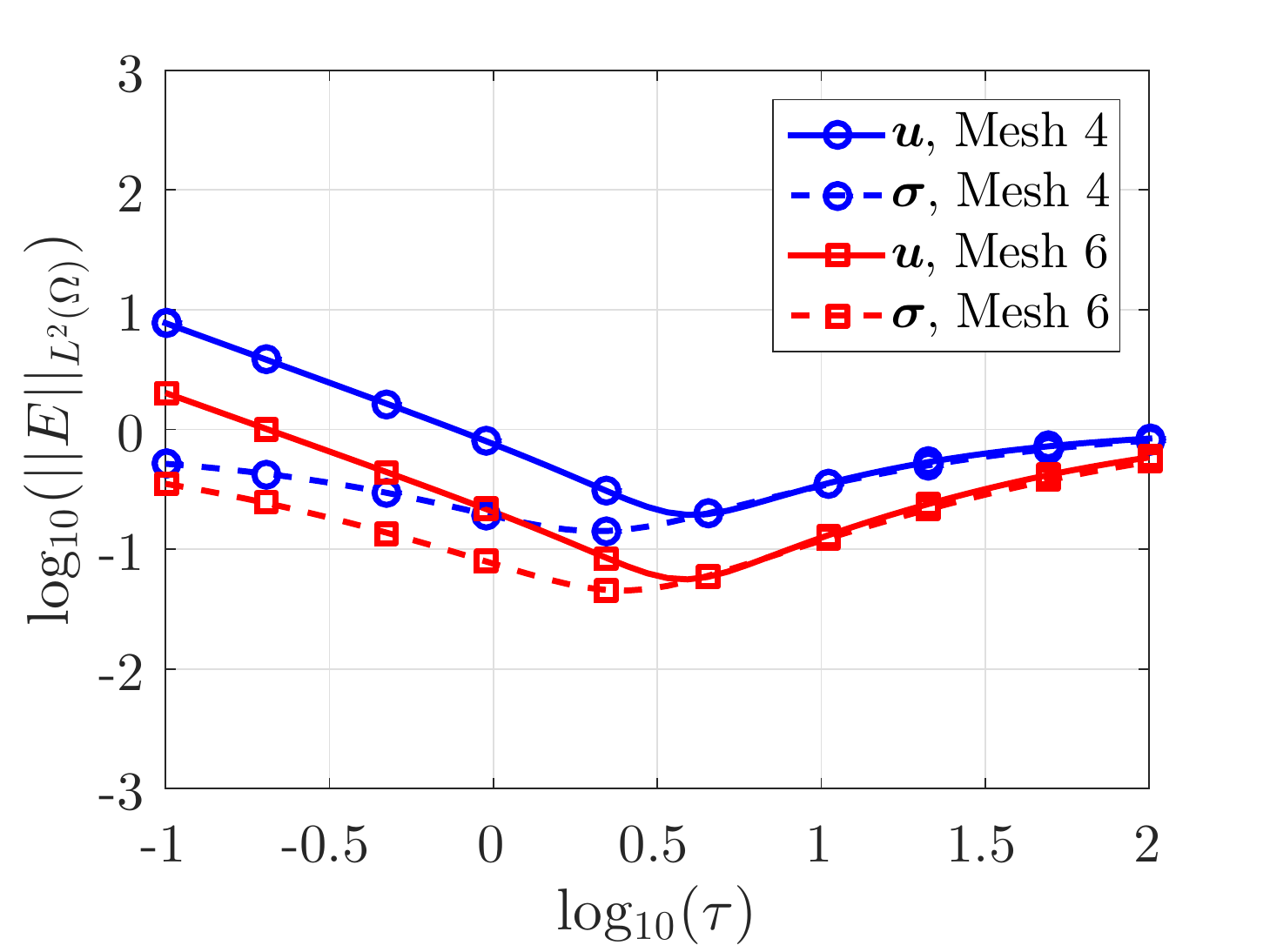}}
	\subfigure[Triangles]{\includegraphics[width=0.49\textwidth]{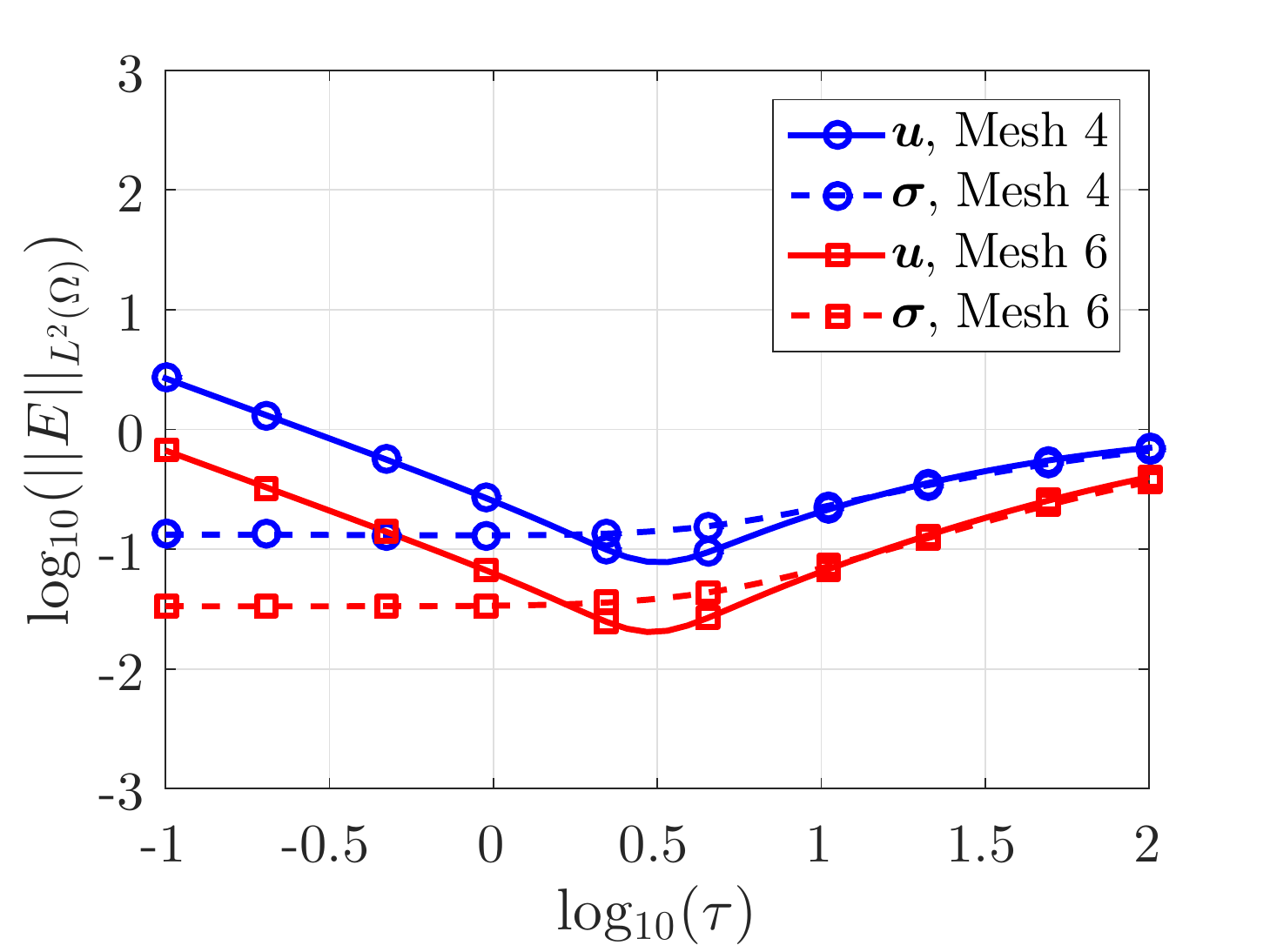}}	
	\caption{Error of the displacement and the stress in the $\eltwo(\Omega)$ norm as a function of the stabilisation parameter $\tau$.}
	\label{fig:tauInfluence2D}
\end{figure}

The proposed methodology requires the choice of the stabilisation tensor $\btau$. Previous works by Cockburn and co-workers~\cite{Jay-CGL:09,Cockburn-CDG:08,soon2009hybridizable} have shown that the stabilisation can have a sizeable effect on the accuracy, convergence and stability of the HDG method. To illustrate the effect of this parameter, the stabilisation tensor is selected as $\btau = \tau (E/\ell) \Insd$, where $\ell$ is a characteristic length. The evolution of the error of the displacement and the stress in the $\eltwo(\Omega)$ norm as a function of the stabilisation parameter $\tau$ is represented in Figure~\ref{fig:tauInfluence2D} for two different computational meshes and for both quadrilateral and triangular elements.
It can be observed that there is an optimum value of the stabilisation parameter, approximately $\tau=3$. It is worth noting that the optimum value is independent on the discretisation considered as the same value provides the most accurate results for all levels of mesh refinement and for all types of element. In addition, the value obtained here for the linear elastic problem also coincides with the optimal value reported in~\cite{doi:10.1002/nme.5833} for the solution of heat transfer problems.

%==========================================================================
\subsection{Influence of the mesh distortion}
\label{sc:distrotion2D}
%==========================================================================

Traditional finite volume methods (e.g. cell-centred and vertex-centred) are well known to suffer an important loss of accuracy when the mesh contains highly distorted elements~\cite{diskin2008,diskin2010comparison,diskin2011comparison,diskin2012}. The accuracy of the reconstruction of the displacement field, required to produce an accurate stress field, is severely compromised by the presence of low quality elements. 

\begin{figure}[!tb]
	\centering
	\subfigure[Mesh 3]{\includegraphics[width=0.24\textwidth]{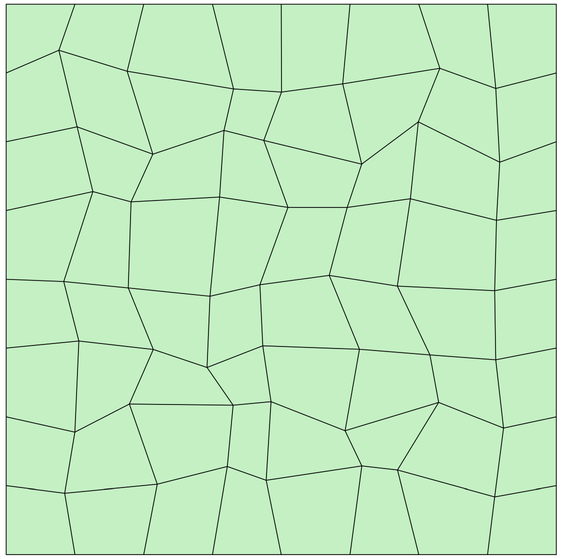}}
	\subfigure[Mesh 5]{\includegraphics[width=0.24\textwidth]{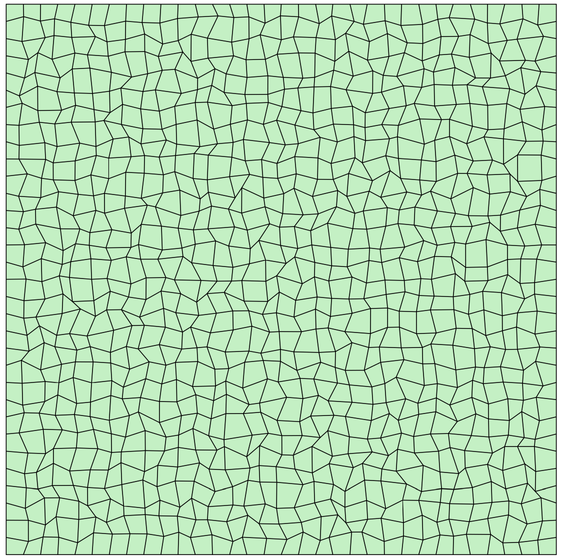}}
	\subfigure[Mesh 3]{\includegraphics[width=0.24\textwidth]{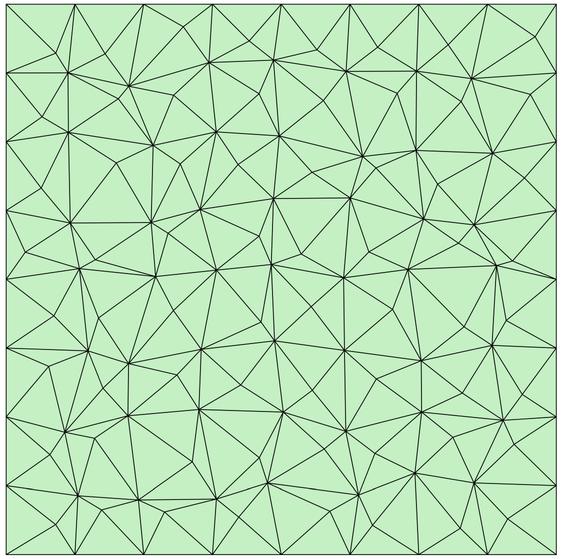}}
	\subfigure[Mesh 5]{\includegraphics[width=0.24\textwidth]{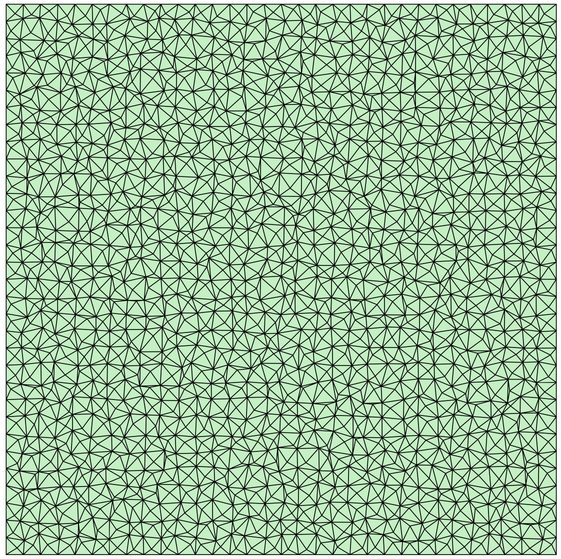}}
	\caption{Distorted quadrilateral and triangular meshes of the domain $\Omega=[0,1]^2$.}
	\label{fig:meshesIrregular}
\end{figure}

To illustrate the robustness of the FCFV method in highly distorted meshes, a new set of meshes is produced by introducing a perturbation of the interior nodes of the uniform meshes employed in the previous example. The new position of an interior node is computed as $\tilde{\bx}_i = \bx_i + \bm{r}_i$, where $\bx_i$ denotes the position in the original uniform mesh and $\bm{r}_i \in \mathbb{R}^{\nsd}$ is a vector containing random numbers generated within the interval $[-h_{\text{min}}/3,h_{\text{min}}/3]$, with $h_{\text{min}}$ being the minimum edge length of the uniform mesh. Two of the meshes with highly distorted elements produced with this strategy are shown in Figure~\ref{fig:meshesIrregular}, for both quadrilateral and triangular elements.	

\begin{figure}[!tb]
	\centering
	\subfigure[Displacement]{\includegraphics[width=0.49\textwidth]{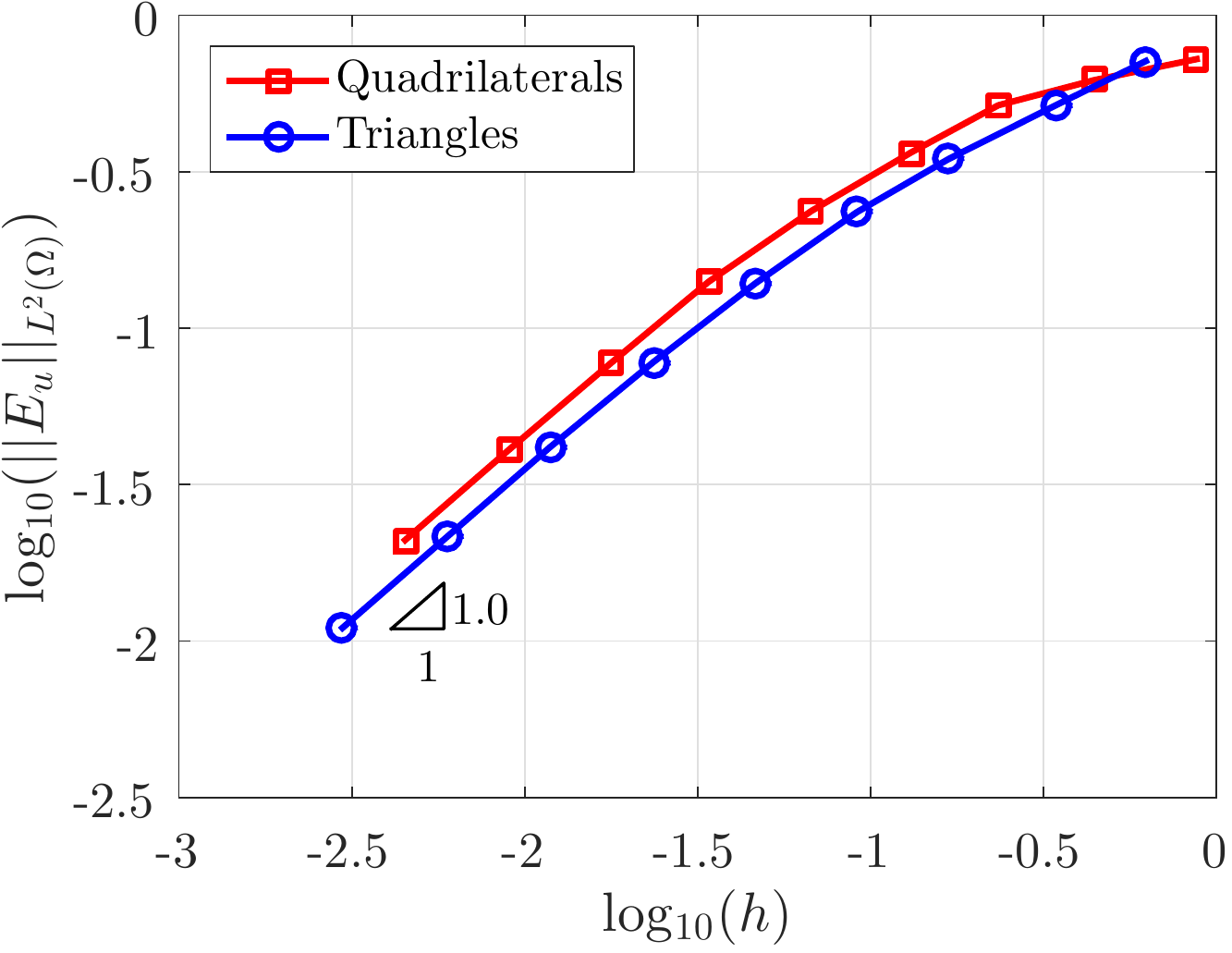}}
	\subfigure[Stress]{\includegraphics[width=0.49\textwidth]{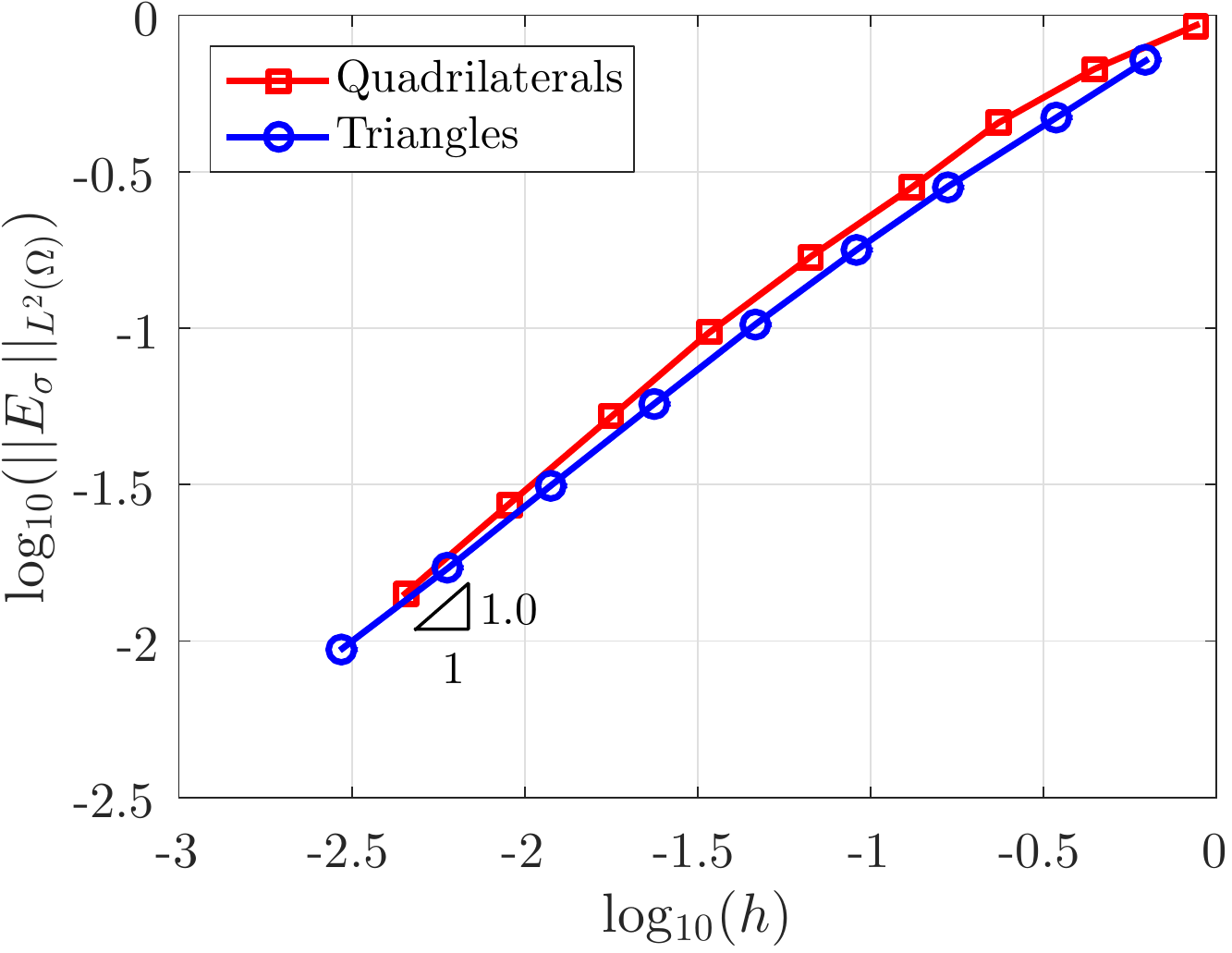}}	
	\caption{Mesh convergence of the $\eltwo(\Omega)$ error of the (a) displacement and (b) the stress, for highly distorted quadrilateral and triangular elements in a medium with $\nu=0.49999$.}
	\label{fig:convergence2DDistored_nu0p49999}
\end{figure}

The numerical experiment of Section~\ref{sc:incompressible} (i.e. for a nearly incompressible medium) is repeated using the distorted quadrilateral and triangular meshes.  Figure~\ref{fig:convergence2DDistored_nu0p49999} shows the error of the computed displacement and stress fields in the $\eltwo(\Omega)$ norm  as a function of the characteristic element size, computed as the maximum of the element diameters in the mesh. 
The results show the expected optimal rate of convergence for both the displacement and the stress, clearly demonstrating that the convergence properties of the proposed approach do not depend upon the quality of the mesh. In addition, by comparing the results of Figures~\ref{fig:convergence2DDistored_nu0p49999} and \ref{fig:convergence2D_nu0p49999}, it can be concluded that the large distortion introduced in the mesh does not result in a sizeable loss of accuracy. For instance, the $\eltwo(\Omega)$ error of the displacement on the uniform mesh of triangular elements in the three finest meshes is 0.037, 0.019 and 0.010 respectively, whereas the error on the corresponding distorted meshes is 0.042, 0.022 and 0.011 respectively.

%==========================================================================
\subsection{Kirch's plate problem}
\label{sc:plate}
%==========================================================================

The next example considers the computation of the stress field in an infinite plate with a circular hole subject to a uniform tension of magnitude $\sigma_0$ in the horizontal direction, a classical test case for solid mechanics solvers in two dimensions~\cite{timoshenko1970theory,szabo1991finite}. The exact solution of the problem is given in polar coordinates by
\begin{subequations}
	\begin{align}
	u_1(r,\theta) & = \frac{\sigma_0 a}{ 8\mu} \left\{ (k+1)\frac{r}{a} \left( 1 + \frac{2a^2}{r^2} \right) \cos(\theta) +  \frac{2a}{r} \left(1 - \frac{a^2}{r^2} \right) \cos(3 \theta) \right\} \\
	u_2(r,\theta) & = \frac{\sigma_0 a}{ 8\mu} \left\{ \frac{r}{a} \left( (k-3) - (k-1)\frac{2a^2}{r^2} \right) \sin(\theta) +  \frac{2a}{r} \left(1 - \frac{a^2}{r^2} \right) \sin(3 \theta) \right\} 
	\end{align}
\end{subequations}
where $\mu$ is the shear modulus and the Kolosov's constant is defined as $k = (3 - \nu)/(1+\nu)$ for plane stress and $k = 3 - 4\nu$ for plane strain.

\begin{figure}[!tb]
	\centering
	\includegraphics[width=0.54\textwidth]{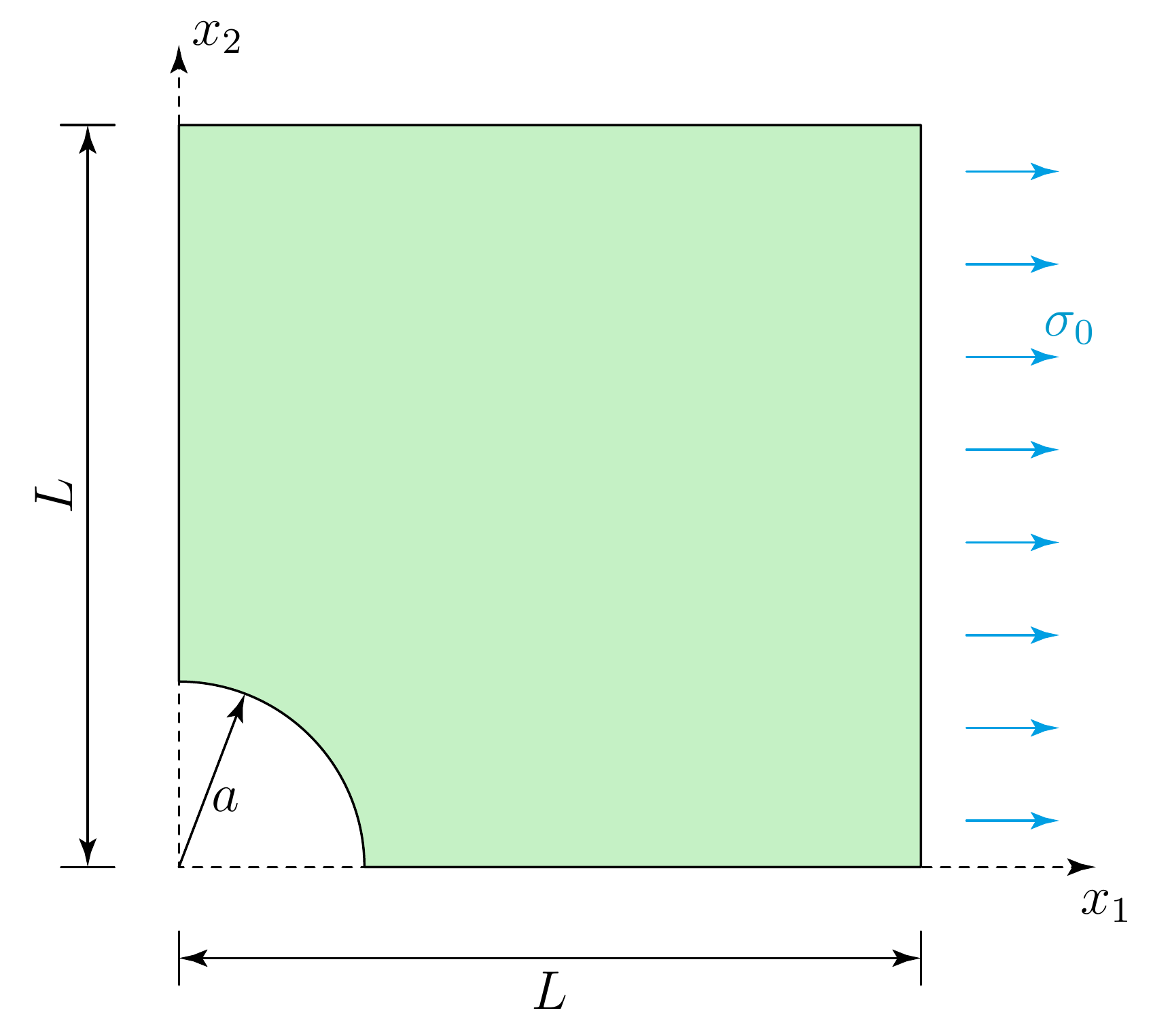}
	\caption{Computational domain for the Kirch's plate problem.}
	\label{fig:KirschPlate}
\end{figure}

The finite computational domain is selected as $[-L,L]^2 \setminus \mathcal{D}_{0,a}$, where $\mathcal{D}_{0,a}$ denotes the disk of radius $a$ centred at the origin. Using the symmetry of the problem, only a quarter of the domain is considered, as illustrated in Figure~\ref{fig:KirschPlate}.

For the numerical examples, $L=4$m, $a=1$m, $E=10^5$Pa, $\nu=0.3$ and $\sigma_0 = 10$Pa are considered. To avoid any effect from the truncation of the infinite domain, the exact traction is imposed on the right and top boundaries, zero traction is imposed on the circular boundary and symmetry boundary conditions are imposed on the bottom and left boundaries.

\begin{figure}[!tb]
	\centering
	\subfigure[Plane stress]{\includegraphics[width=0.49\textwidth]{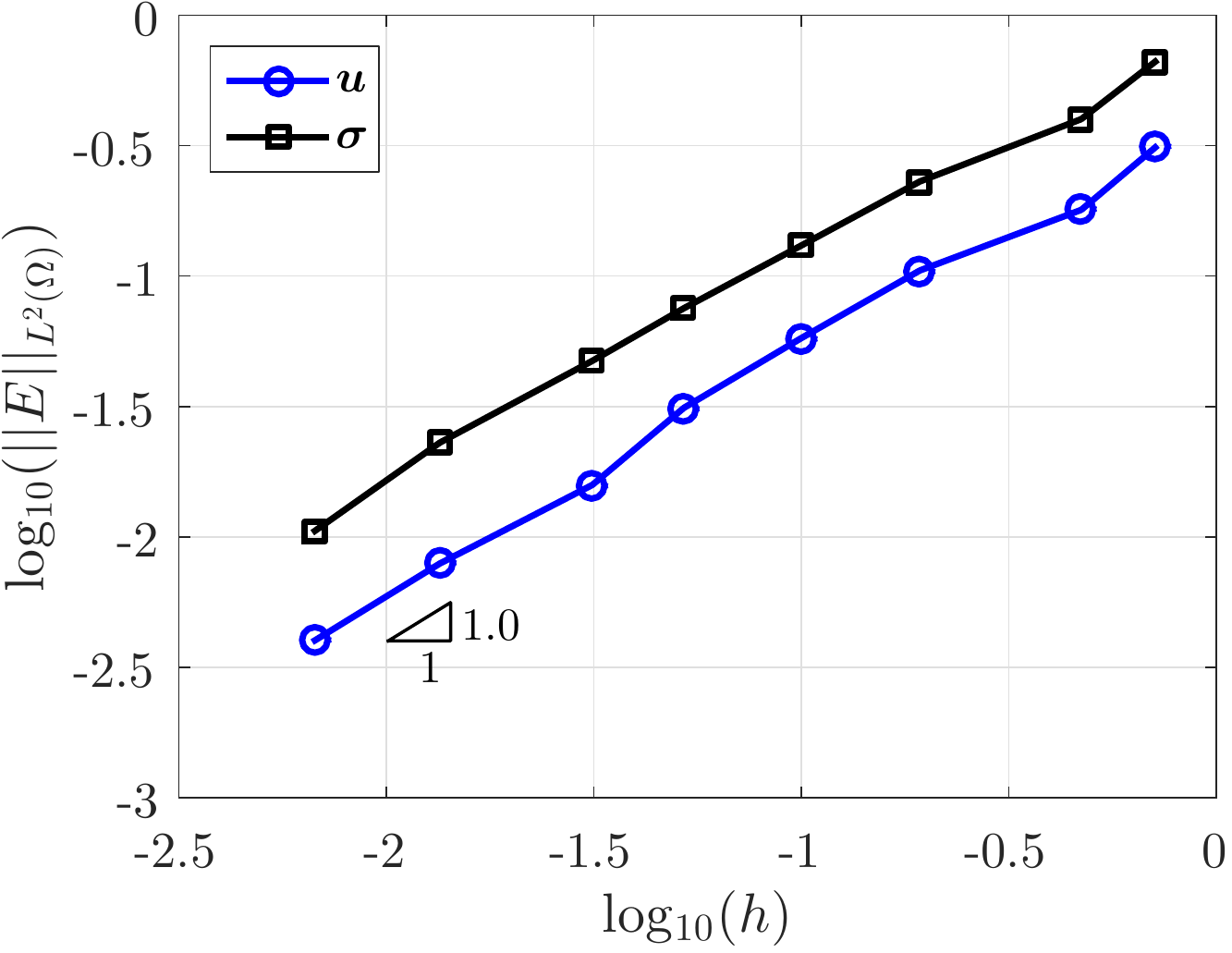}}
	\subfigure[Plane strain]{\includegraphics[width=0.49\textwidth]{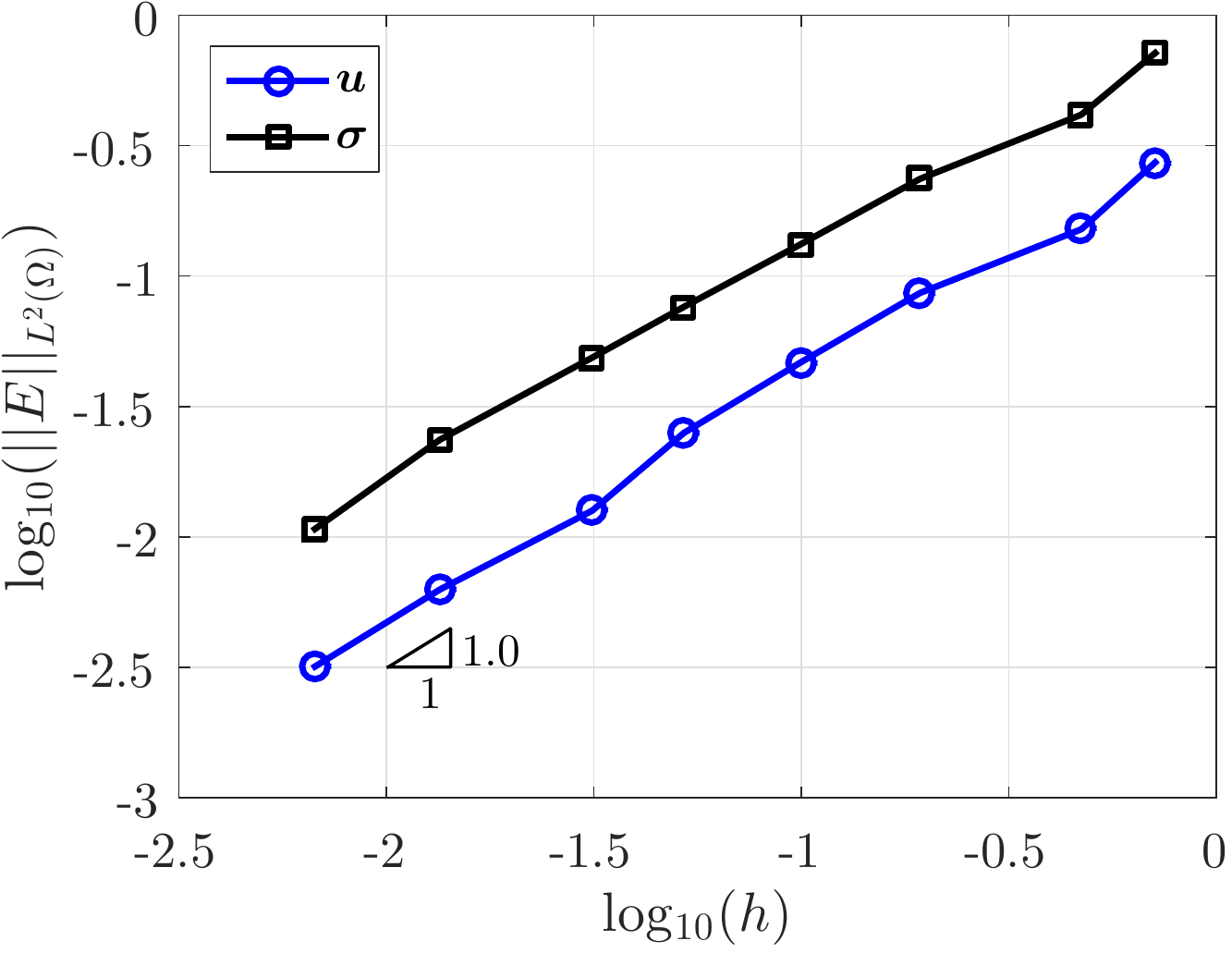}}	
	\caption{Kirch's plate problem: mesh convergence of the $\eltwo(\Omega)$ error of the  displacement and the stress for a (a) plane stress and (b) plane strain two dimensional models.}
	\label{fig:KirschPlateHConv}
\end{figure}

The convergence of the displacement and stress error measured in the $\eltwo(\Omega)$ norm as a function of the characteristic element size is shown in Figure~\ref{fig:KirschPlateHConv}, showing the expected rate for both quantities. 

\begin{figure}[!tb]
	\centering
	\subfigure[$\sigma_{11}$]{\includegraphics[width=0.32\textwidth]{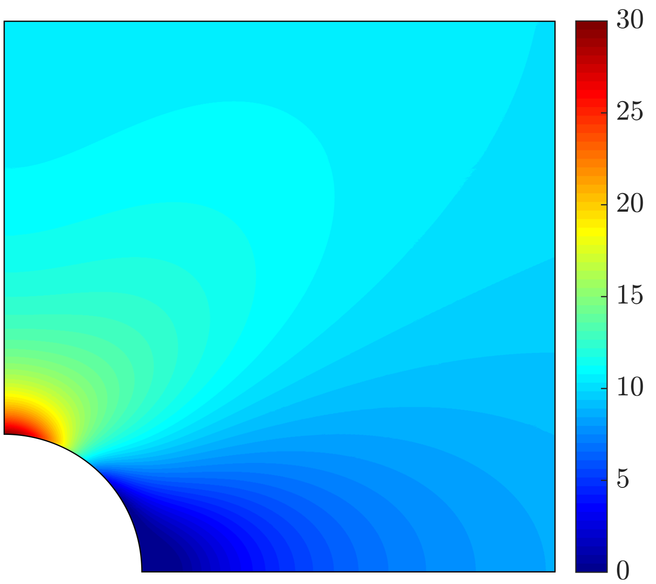}}
	\subfigure[$\sigma_{22}$]{\includegraphics[width=0.32\textwidth]{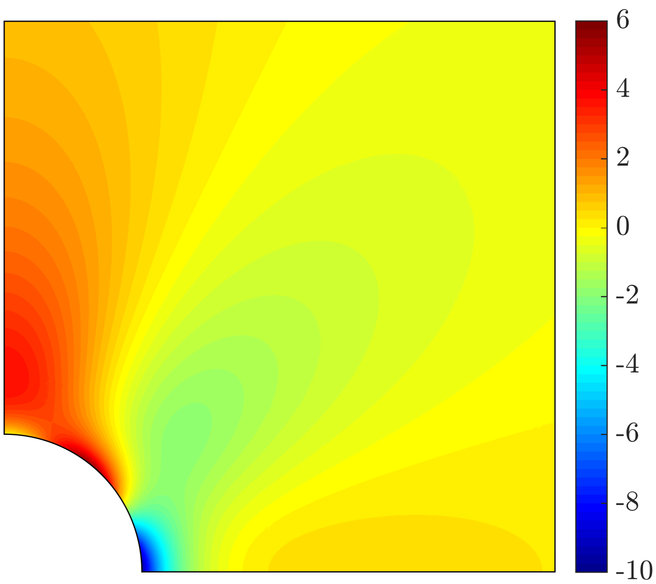}}	
	\subfigure[$\tau_{12}$]{\includegraphics[width=0.32\textwidth]{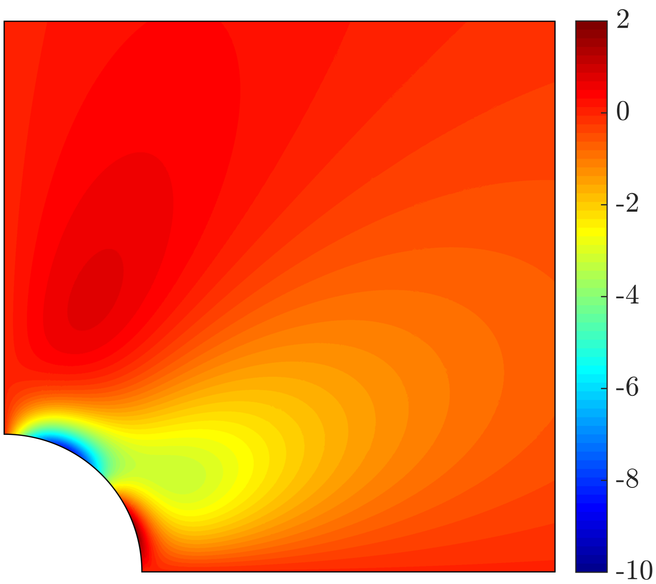}}	
	\caption{Kirch's plate problem: computed components of the stress field.}
	\label{fig:KirschPlateStress}
\end{figure}

The stress field computed on the seventh mesh used for the mesh convergence study, with 573,123 triangular elements, is shown in Figure~\ref{fig:KirschPlateStress}. 
The computation with the proposed FCFV method required the solution of a linear system of 1,719,826 equations, taking 6 seconds to compute all elemental matrices, 3 seconds to perform the assembly of the global system and 46 seconds to solve using a direct method. The developed code is written in Matlab and the computation was performed in an Intel$^{\tiny{\textregistered}}$
Xeon$^{\tiny{\textregistered}}$ CPU $@$ 3.70GHz and 32GB main memory available.

%==========================================================================
\subsection{Cook's membrane problem}
\label{sc:cook}
%==========================================================================

\begin{figure}
	\centering
	\includegraphics[width=0.4\textwidth]{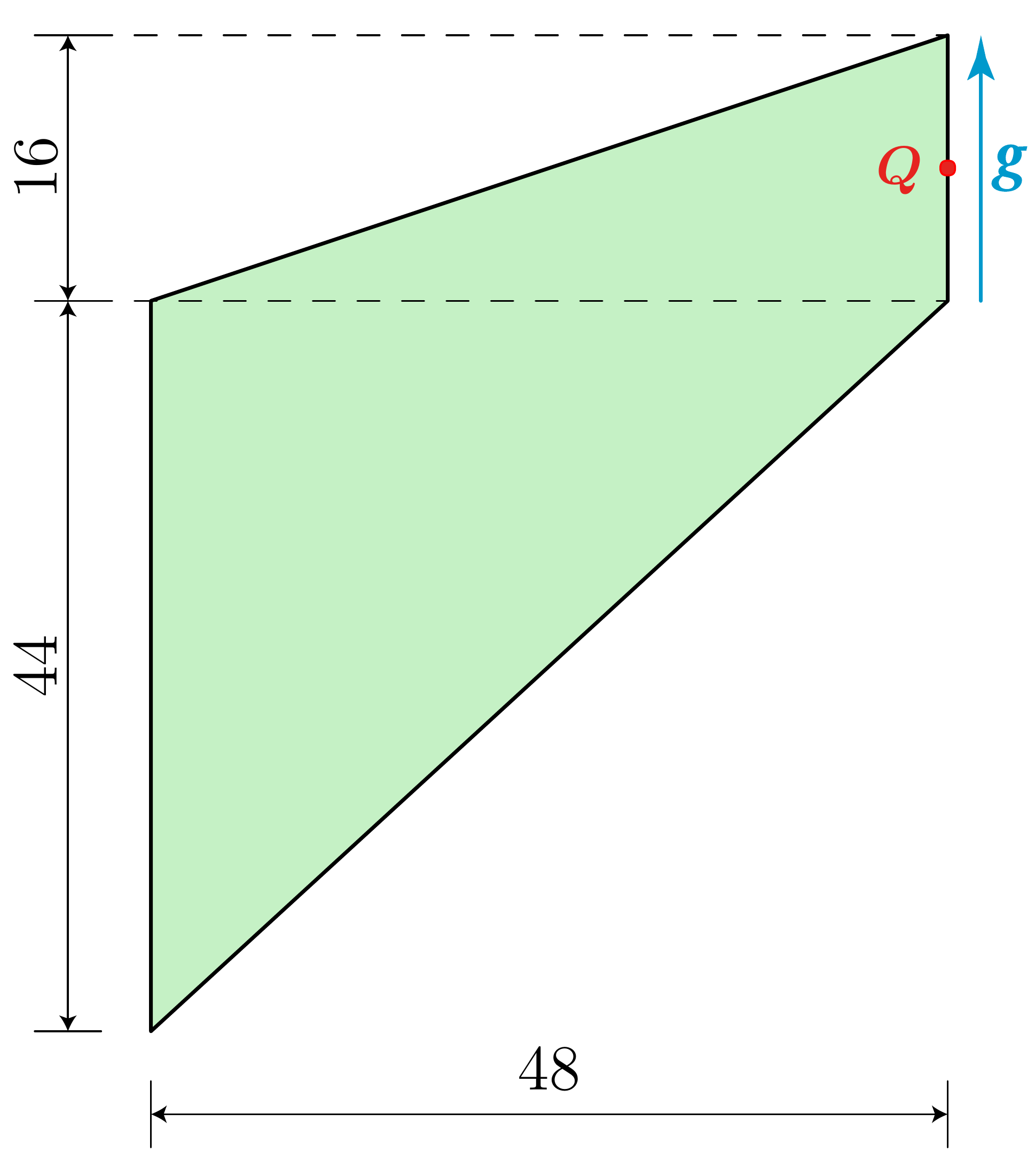}
	\caption{Cook's membrane problem.}
	\label{fig:cooksMembraneSetup}
\end{figure}

The last two dimensional example considers a classical bending dominated test case employed to validate the susceptibility of linear elastic solvers to volumetric locking, the so-called Cook's membrane problem~\cite{cook2001concepts}. The problem consists of a tapered plate clamped on one end and subject to a shear load, taken as $\bg=(0,1/16)$ here, on the opposite end, as illustrated in Figure~\ref{fig:cooksMembraneSetup}. 

Two cases, reported in~\cite{auricchio2005analysis}, are considered to validate the performance of the recently proposed FCFV methodology. The first case involves a material with Young modulus $E=1$ and Poisson ratio $\nu=1/3$ and the second case a nearly incompressible material with Young modulus $E = 1.12499998125$ and Poisson ratio $\nu = 0.499999975$. As there is no analytical solution available, the vertical displacement at the mid point of the right end of the plate, $\bm{Q} =(48,52)$, is compared against the reference values reported in~\cite{auricchio2005analysis}, given by 21.520 and 16.442 respectively.

\begin{figure}[!tb]
	\centering
	\subfigure[$\nu = 1/3$]{\includegraphics[width=0.49\textwidth]{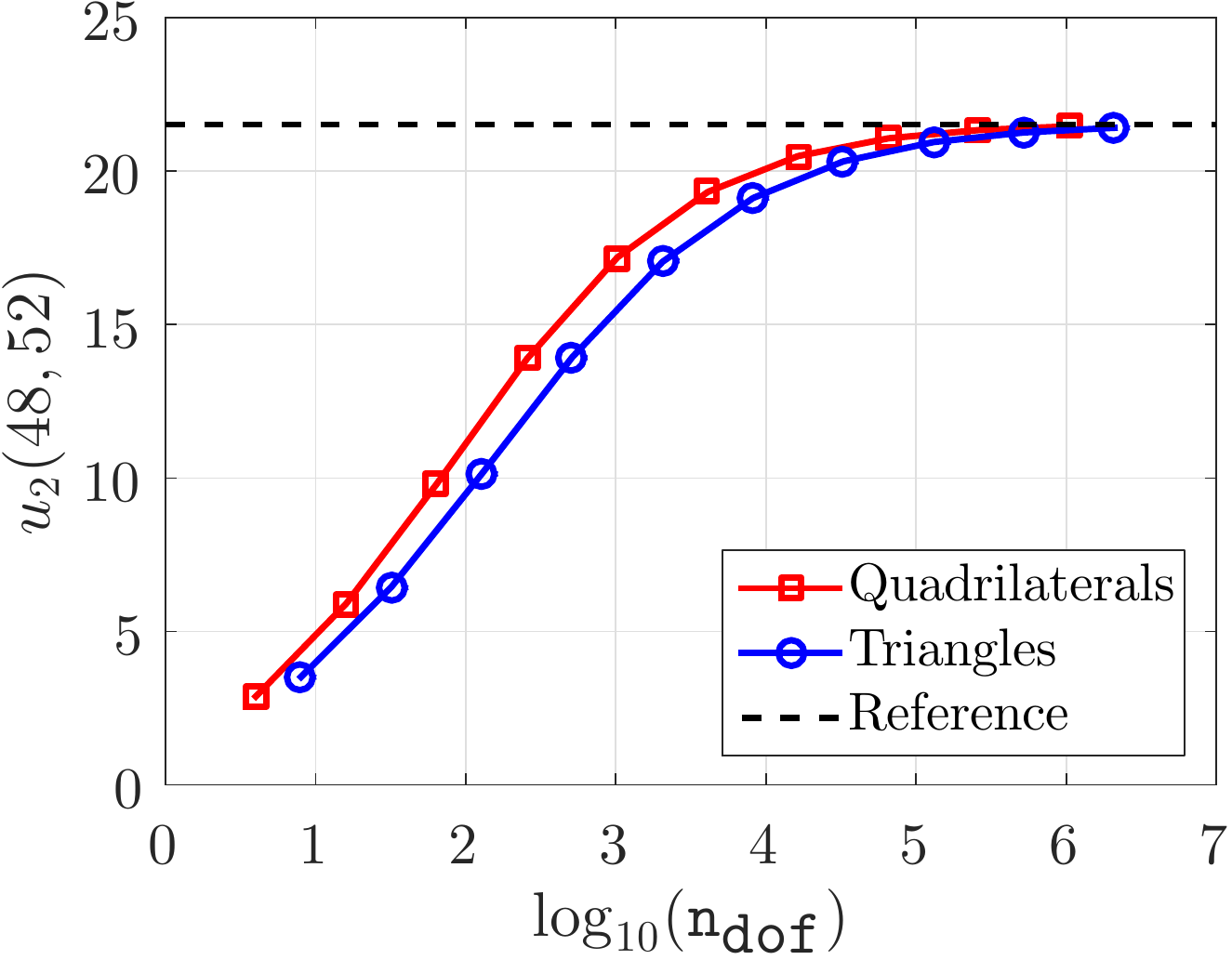}}
	\subfigure[$\nu = 0.499999975$]{\includegraphics[width=0.49\textwidth]{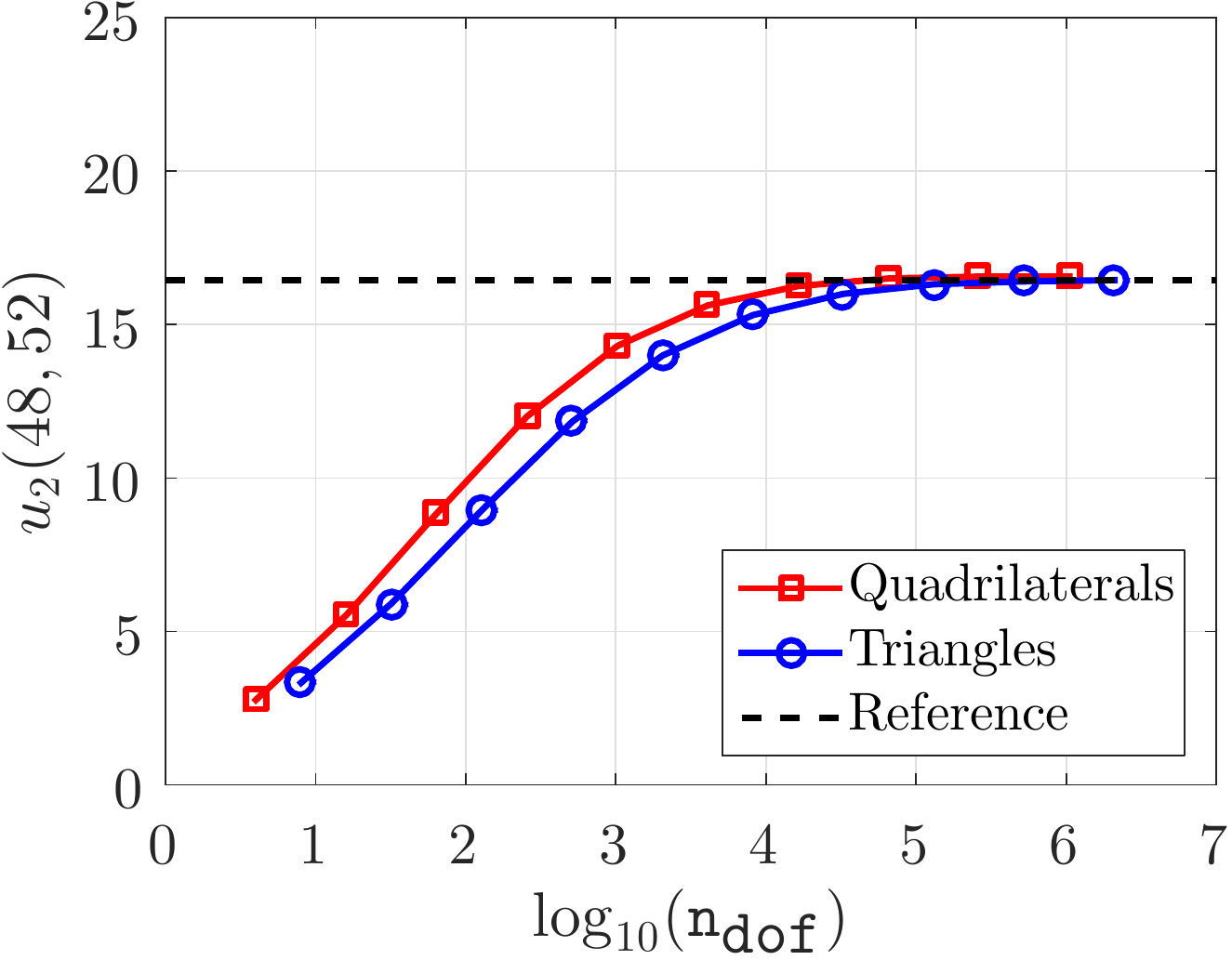}}	
	\caption{Cook's membrane problem: evolution of the vertical displacement at the mid point of the right end of the plate as a function of the total number of degrees of freedom.}
	\label{fig:cooksMembrane}
\end{figure}

\begin{figure}[!tb]
	\centering
	\subfigure[$u_1$]{\includegraphics[width=0.32\textwidth]{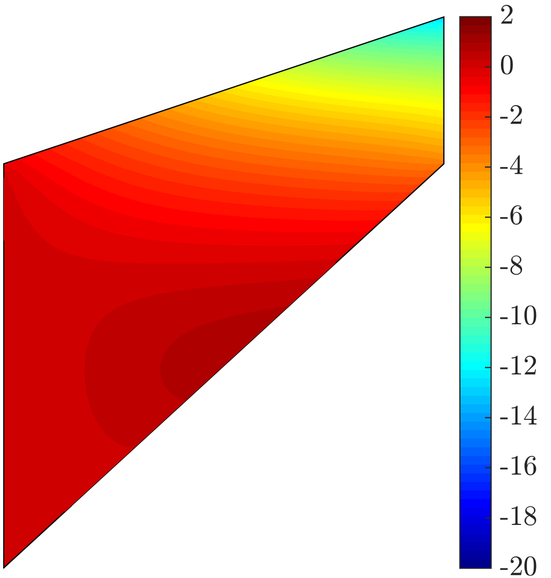}}
	\subfigure[$u_2$]{\includegraphics[width=0.32\textwidth]{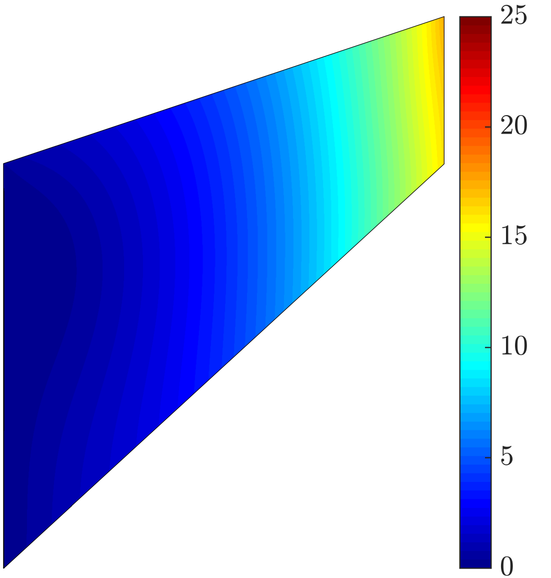}}	
	\subfigure[$\sigma_\text{VM}$]{\includegraphics[width=0.32\textwidth]{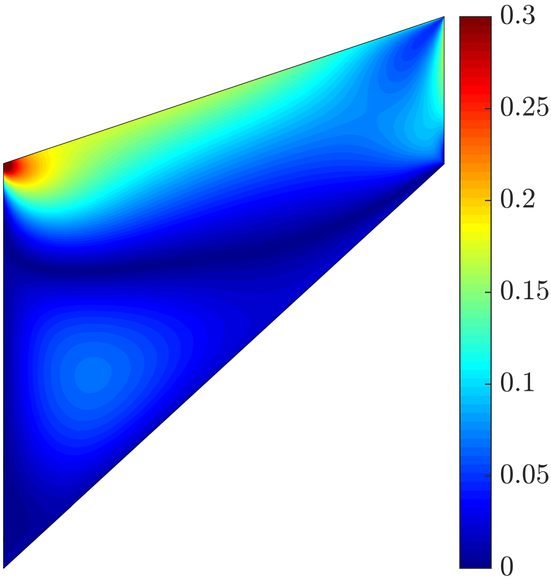}}	
	\caption{Cook's membrane problem: displacement and Von Mises stress.}
	\label{fig:cooksMembraneCase2}
\end{figure}

Figure~\ref{fig:cooksMembrane} shows the convergence of the vertical displacement at point $\bm{Q}$ for both cases and using both quadrilateral and triangular elements. 
The results indicate convergence of the vertical displacement in the ninth mesh, with 262,144 elements, for the first case with $\nu = 1/3$. The computed displacement at the mid point of the right end of the plate is within a 1\% difference with respect to the results reported in~\cite{auricchio2005analysis}. The FCFV computation required the solution of a linear system of 1,049,600 equations, taking 4 seconds to compute all elemental matrices, 2 seconds to perform the assembly of the global system and 1 minute to solve using a direct method. The displacement field and Von Mises stress for this computation are represented in Figure~\ref{fig:cooksMembraneCase2}.

For the second case, with a nearly incompressible material, the results in the eight mesh, with 65,536 elements, show convergence to the reference value, illustrating the robustness and accuracy of the proposed approach in the incompressible limit. The computed displacement at the mid point of the right end of the plate is within a 0.5\% difference with respect to the results reported in~\cite{auricchio2005analysis}. The FCFV computation required the solution of a linear system of 262,656 equations, taking 1 second to compute all elemental matrices, 0.5 seconds to perform the assembly of the global system and 10 seconds to solve using a direct method.

%==========================================================================
\section{Three dimensional examples}
\label{sc:examples32D}
%==========================================================================

This Section presents three numerical examples in three dimensions to show the potential of the proposed FCFV approach in more complicated scenarios, including a more realistic application involving a complex geometry.

%==========================================================================
\subsection{Cantilever beam under shear}
\label{sc:beamBending}
%==========================================================================

The first three-dimensional example involves the analysis of a beam under shear and it is used here to verify the optimal convergence properties of the FCFV method in three dimensions for both hexahedral and tetrahedral elements. 

The analytical solution of the problem is given by~\cite{barber2002elasticity}
\begin{subequations}
	\begin{align}
	u_1(x_1,x_2,x_3) & = \frac{3P\nu}{4E} x_1 x_2 x_3,\\
	u_2(x_1,x_2,x_3) & = \frac{P}{8E} \left[ 3 \nu x_3 \left( x_1^2 - x_2^2 \right) - x_3^3 \right], \\
	\begin{split}
	u_3(x_1,x_2,x_3) & = \frac{Px_2}{8E} \left[  \nu(3x_1^2 - x_2^2 + 4) + 3x_3^2 - 2x_2^2 + 6   \right] \\ & - \frac{3P\nu}{\pi^3 E} \sum_{n=1}^\infty \frac{(-1)^n}{n^3 \cosh(n\pi)}\cos(n\pi x)\sinh(n\pi y).
	\end{split}
	\end{align}
\end{subequations}

The domain is $\Omega=[-1,1] \times [-1,1] \times [0,L]$ and the material properties are taken as $E=25$ and $\nu=0.3$. Following~\cite{gain2014virtual}, the boundary conditions correspond to the exact displacement imposed on $\Gamma_D = \{(x_1,x_2,x_3) \in \mathbb{R}^3 \; | \; x_3=L\}$ whereas the exact tractions are imposed at $\Gamma_N = \partial \Omega \setminus \Gamma_D$. The length of beam is $L=10$ and the shear load is taken as $P=0.1$.

\begin{figure}[!tb]
	\centering
	\subfigure[Displacement]{\includegraphics[width=0.49\textwidth]{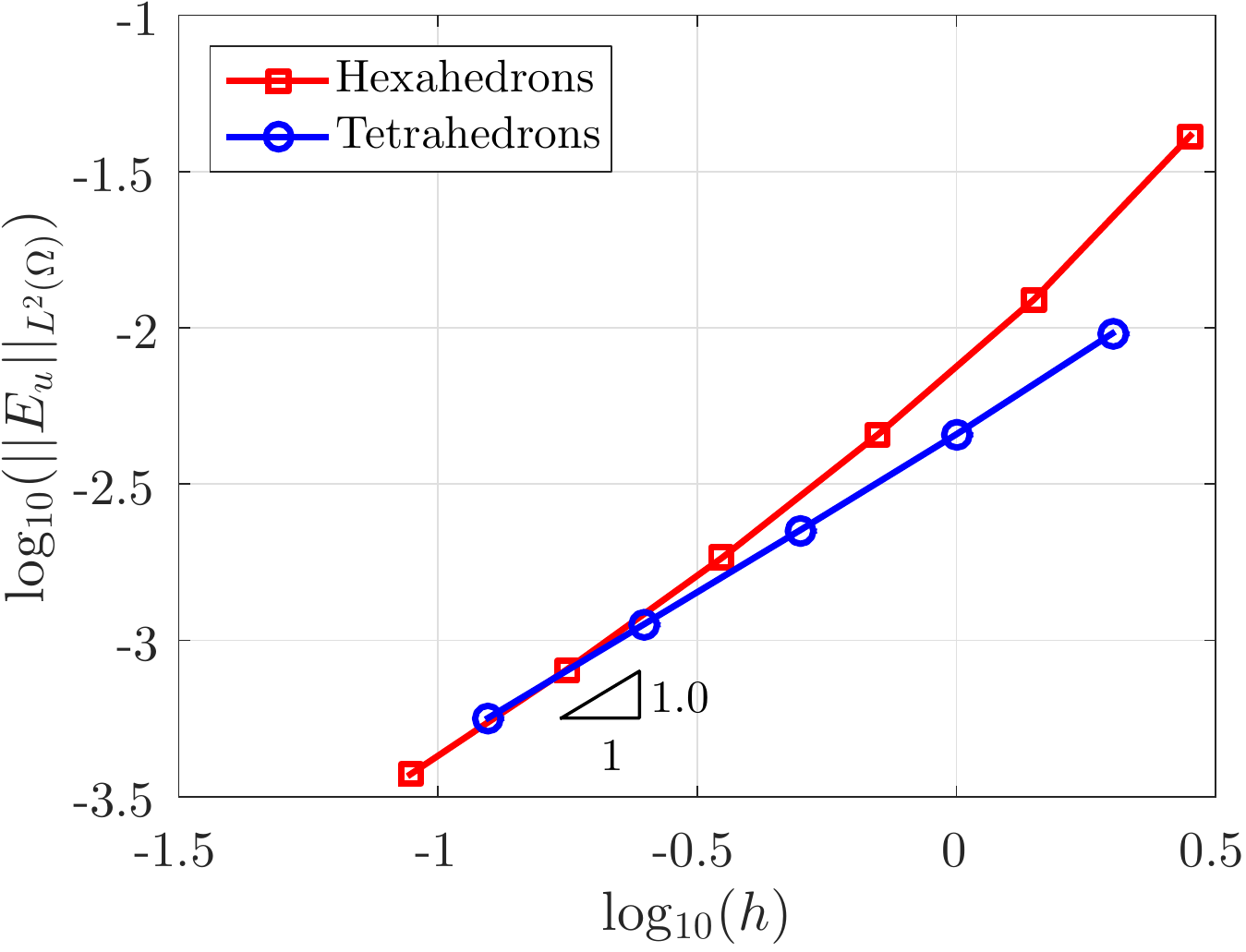}}
	\subfigure[Stress]{\includegraphics[width=0.49\textwidth]{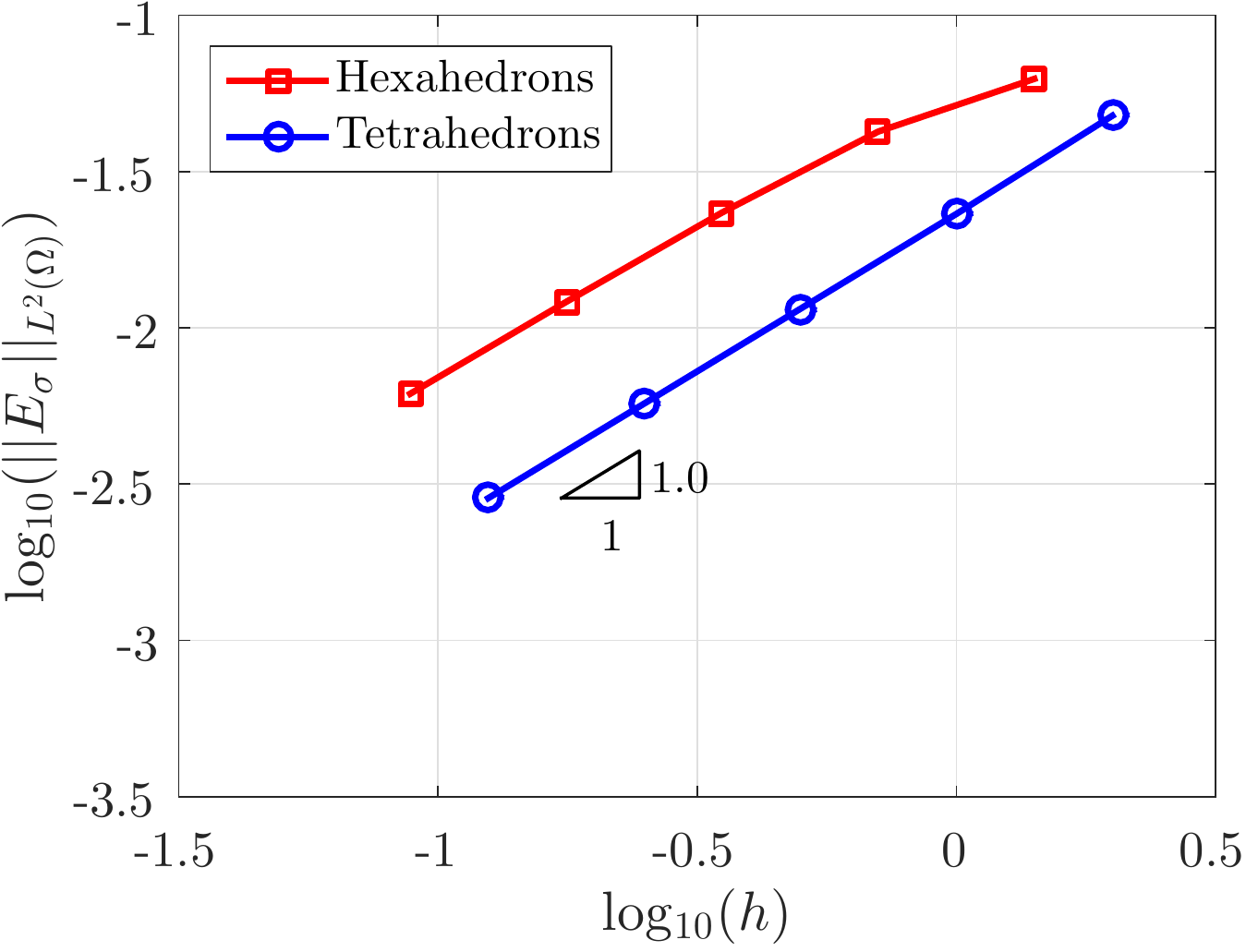}}	
	\caption{Cantilever beam under shear: mesh convergence of the $\eltwo(\Omega)$ error of the (a) displacement and (b) the stress, for hexahedral and tetrahedral elements.}
	\label{fig:beamConvergence}
\end{figure}

Five tetrahedral and six hexahedral meshes are considered to perform the mesh convergence analysis. The tetrahedral meshes contain 120, 960, 7,680, 61,440 and 491,520 elements respectively, whereas the hexahedral meshes contain 5, 40, 320, 2,560, 20,480 and 163,840 elements respectively. Figure~\ref{fig:beamConvergence} displays the error of the computed displacement and the stress fields in the $\eltwo(\Omega)$ norm as a function of the characteristic element size, showing the optimal approximation properties of the proposed FCFV method in three dimensions for both hexahedral and tetrahedral elements. 

\begin{figure}[!tb]
	\centering
	\subfigure[$u_1$]{\includegraphics[width=0.24\textwidth]{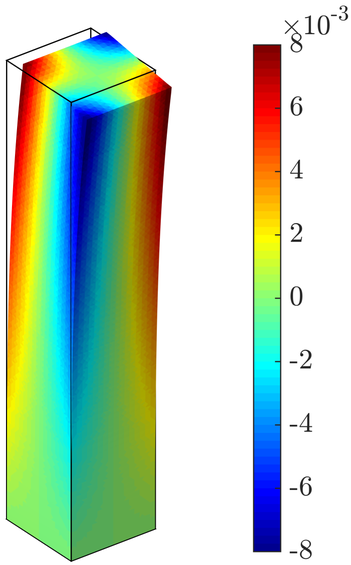}}
	\subfigure[$u_2$]{\includegraphics[width=0.24\textwidth]{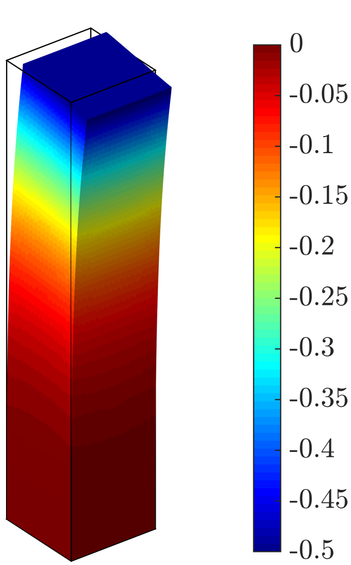}}	
	\subfigure[$u_3$]{\includegraphics[width=0.24\textwidth]{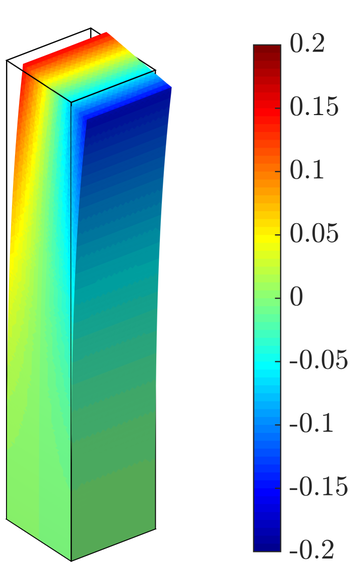}}	
	\subfigure[$\sigma_\text{VM}$]{\includegraphics[width=0.24\textwidth]{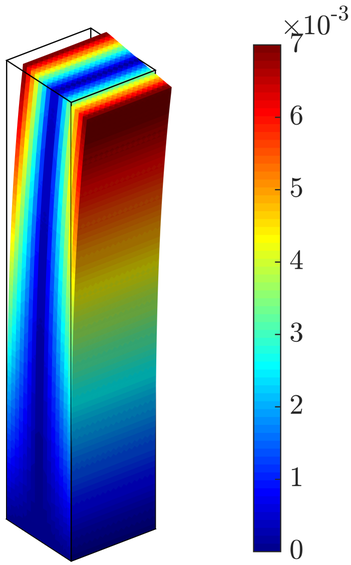}}	
	\caption{Cantilever beam under shear: displacement and Von Mises stress.}
	\label{fig:beamSols}
\end{figure}

The three components of the displacement and the Von Mises stress are represented in Figure~\ref{fig:beamSols}. 
The results, corresponding to the finer tetrahedral mesh are displayed over the deformed configuration. The computation with the proposed FCFV method required the solution of a linear system of 2,915,328 equations, taking 10 seconds to compute all elemental matrices, 5 seconds to perform the assembly of the global system and 10 minutes to solve using a direct method. The developed code is written in Matlab and the computation was performed in an Intel$^{\tiny{\textregistered}}$ Xeon$^{\tiny{\textregistered}}$ CPU $@$ 3.70GHz and 32GB main memory available.

%==========================================================================
\subsection{Thin cylindrical shell}
\label{sc:cylindricalShell}
%==========================================================================

The next example involves the analysis of a thin cylindrical shell subject to a uniform internal pressure and with fixed ends. This is a particularly challenging problem for low order methods due to the localised bending occurring near the ends of the shell, leading to a radial displacement that exhibits a boundary layer behaviour.  

The analytical solution of the problem is given by~\cite{timoshenko1959theory}
\begin{subequations}
	\begin{align}
	u_1(x_1,x_2,x_3) & = u_r(x_3), \cos(\theta) \\
	u_2(x_1,x_2,x_3) & = u_r(x_3), \sin(\theta) \\
	u_3(x_1,x_2,x_3) & = 0,
	\end{align}
\end{subequations}
where $u_r$ denotes the radial displacement, given by
\begin{equation} 
u_r(x_3) = -\frac{P a^2}{E t} \left(1 - C_1 \sin(\beta x_3) \sinh(\beta x_3) - C_2 \cos(\beta x_3) \cosh(\beta x_3) \right),
\end{equation}
with 
\begin{equation}
C_1 =  \frac{2 \sin(\alpha) \sinh(\alpha) }{\cos(2\alpha) + \cosh(2\alpha)},  \qquad 
C_2 =  \frac{2 \cos(\alpha) \cosh(\alpha) }{\cos(2\alpha) + \cosh(2\alpha)}.
\end{equation}
In the above expressions, $P$ is the magnitude of the internal pressure, $a$ is the midplane radius of the shell, $t$ is the thickness, 
\begin{equation}
\alpha = \frac{\beta L}{2}, \qquad \beta = \left( \frac{E t}{4a^2 D} \right)^{1/4}, \qquad D = \frac{Et^3}{12(1-\nu^3)}
\end{equation}
and $L$ is the height of the shell.

\begin{figure}[!tb]
	\centering
	\subfigure[]{\includegraphics[width=0.15\textwidth]{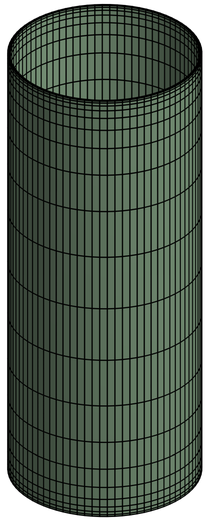}}
	\hspace{2cm}
	\subfigure[]{\includegraphics[width=0.48\textwidth]{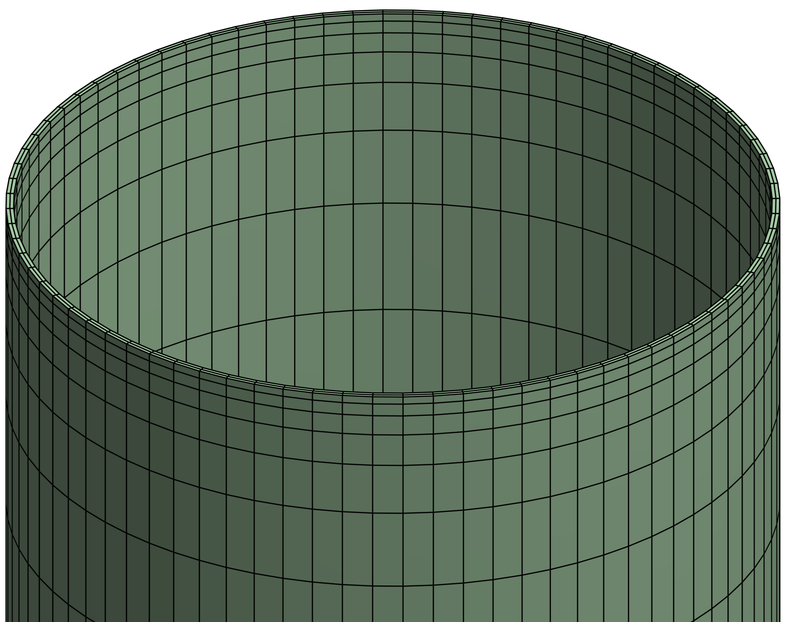}}
	\caption{Thin cylindrical shell: (a) hexahedral mesh with element stretching near the ends of the shell and (b) detail of the mesh showing the two elements across the thickness.}
	\label{fig:shell_Mesh}
\end{figure}

The numerical results presented here consider $L=5$, $a=1$, $t=0.02$, $E=1$ and $\nu=0.3$. Hexahedral meshes with element stretching are considered to capture the localised variation of the displacement near the ends of the shell. Figure~\ref{fig:shell_Mesh} shows one hexahedral mesh with 3,200 elements and a detail of the mesh near the end, illustrating the stretching used and showing that only two elements are considered across the thickness. It is worth emphasising that the problem is solved using solid elements despite the shell theory is applicable in this problem~\cite{timoshenko1959theory}.

\begin{figure}[!tb]
	\centering
	\subfigure[$u_1$]{\includegraphics[width=0.24\textwidth]{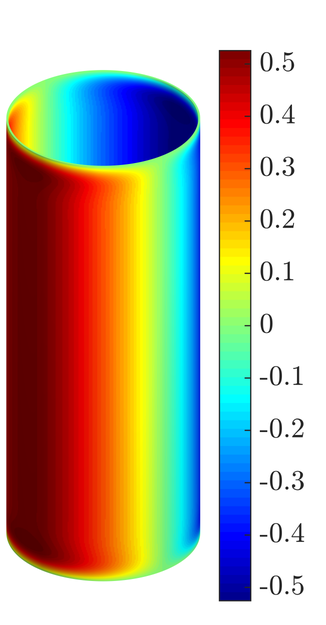}}
	\subfigure[$u_2$]{\includegraphics[width=0.24\textwidth]{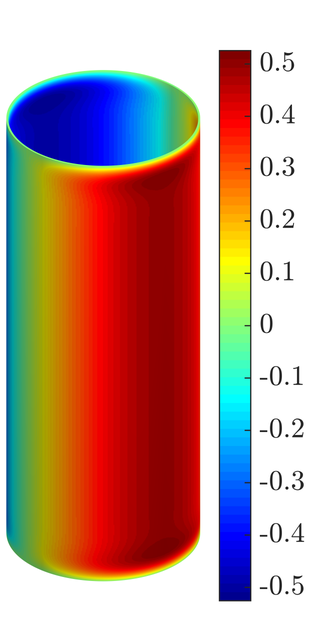}}	
	\subfigure[$u_3$]{\includegraphics[width=0.24\textwidth]{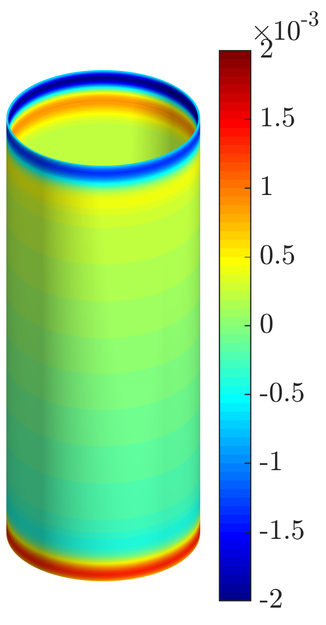}}	
	\subfigure[$u_r$]{\includegraphics[width=0.24\textwidth]{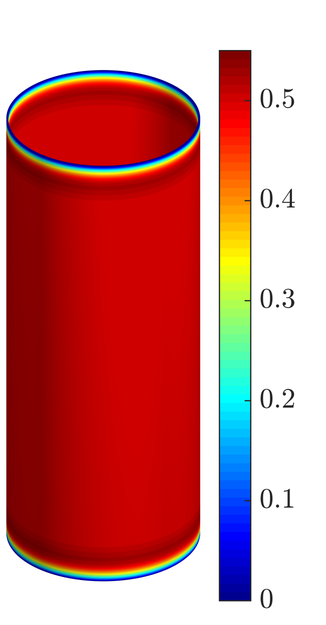}}	
	\caption{Thin cylindrical shell: components of the displacement field and radial displacement.}
	\label{fig:shell_u}
\end{figure}

The three components of the displacement field and the radial displacement are depicted in Figure~\ref{fig:shell_u} on a fine mesh with 819,200 elements.
The results are in excellent agreement with the shell theory, with an $\eltwo(\Omega)$ error of  $9.1 \times 10^{-4}$ and $4.1 \times 10^{-3}$ in the displacement and stress respectively.

\begin{figure}[!tb]
	\centering
	\includegraphics[width=0.5\textwidth]{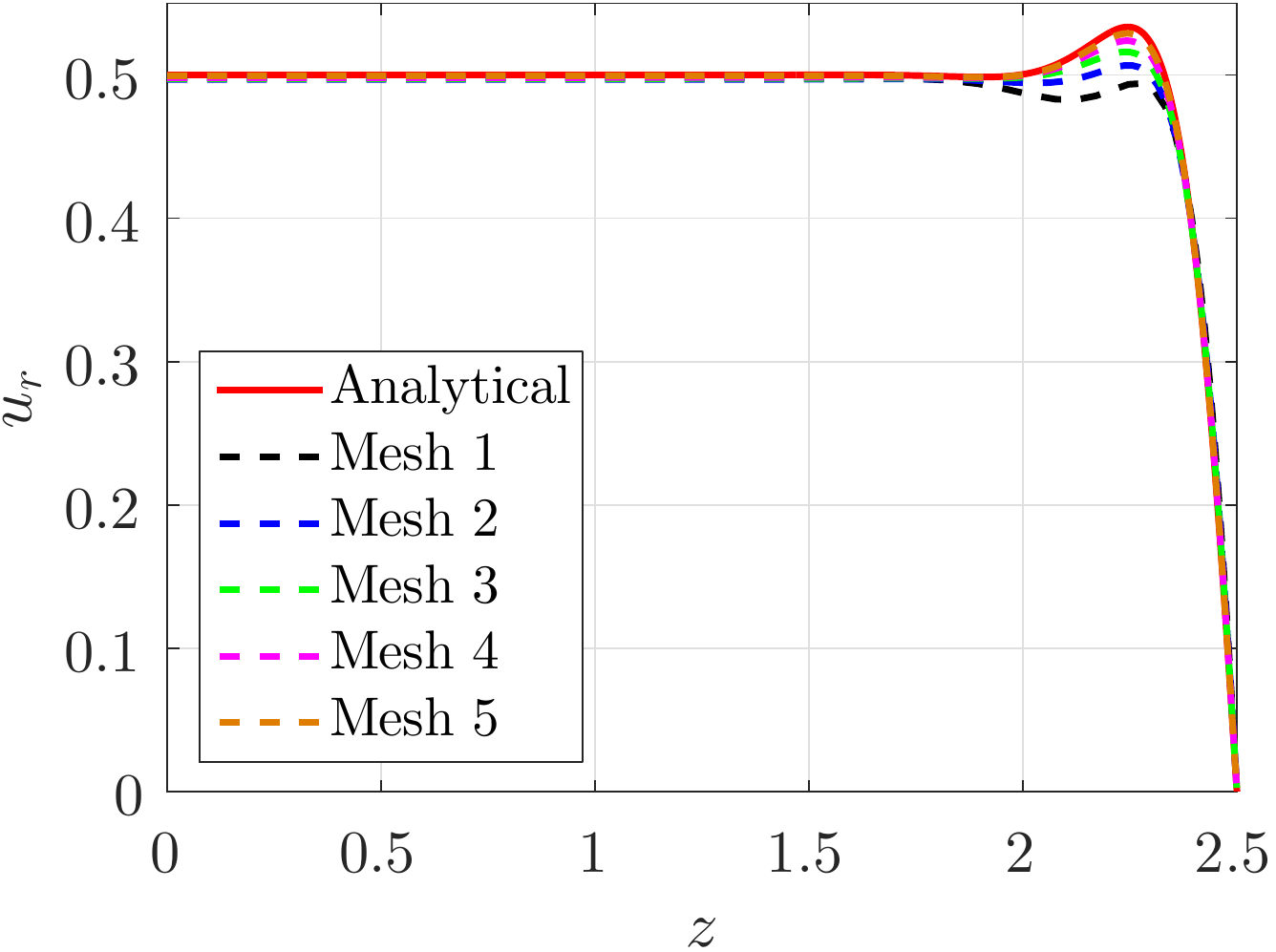}
	\caption{Thin cylindrical shell: mesh convergence of the radial displacement.}
	\label{fig:shell_RadialUsection}
\end{figure}

Figure~\ref{fig:shell_RadialUsection} shows a detailed view of the radial displacement computed with five subsequently refined meshes and compared to the analytical solution.
The results show the ability of the proposed FCFV methodology to capture the boundary layer behaviour of the radial displacement with a mesh with only two elements across the thickness. 

A more quantitative analysis is presented in Table~\ref{tab:shell}, where the number of elements the number of degrees of freedom and the error of the displacement field, the stress field and the radial displacement is given for the five meshes utilised.
\begin{table}[hbt]
	\centering
	\begin{tabular}[hbt]{|c|c|c|c|c|c|}
		\hline
		 Mesh & Elements & $\ndof$ & $E_u$ & $E_{\stress}$ & $E_r$   \\
		\hline
		1 & $80 \times 10 \times 2$    &    33,480 & 0.0055 & 0.0764 & 0.0285 \\
		\hline
		2 & $160 \times 20 \times 2$   &   134,160 & 0.0035 & 0.0409 & 0.0191  \\
		\hline
		3 & $320 \times 40 \times 2$   &   997,440 & 0.0021 & 0.0224 & 0.0124  \\
		\hline
		4 & $640 \times 80 \times 2$   & 3,991,680 & 0.0011 & 0.0116 & 0.0070  \\
		\hline
		5 & $1280 \times 160 \times 2$ & 8,599,680 & 0.0005 & 0.0055 & 0.0031  \\
		\hline
	\end{tabular}
	\caption{Thin cylindrical shell: details of the mesh convergence analysis. For each mesh, the number of elements, the number of degrees of freedom and the error of the displacement field, the stress field and the radial displacement are given.}
	\label{tab:shell}
\end{table}
The error of the radial displacement is measured over a section, corresponding to $x_1=x_1^*$ and $x_2=x_2^*$, as
\begin{equation}
E_r = \left\{ \frac{\int_{-L/2}^{L/2} \left[u_r^h(x_1^*,x_2^*,x_3) - u_r(x_3) \right]^2 dx_3 }{ \int_{-L/2}^{L/2} u_r(x_3)^2 dx_3 } \right\}^{1/2},
\end{equation}
where $u_r^h$ and $u_r$ are the computed and exact radial displacement respectively.

The results in Table~\ref{tab:shell} show, once more, the optimal first-order convergence of the error of the displacement and stress fields under mesh refinement.

%==========================================================================
\subsection{Bearing cap}
\label{sc:bearing}
%==========================================================================

\begin{figure}[!tb]
	\centering
	\includegraphics[width=0.8\textwidth]{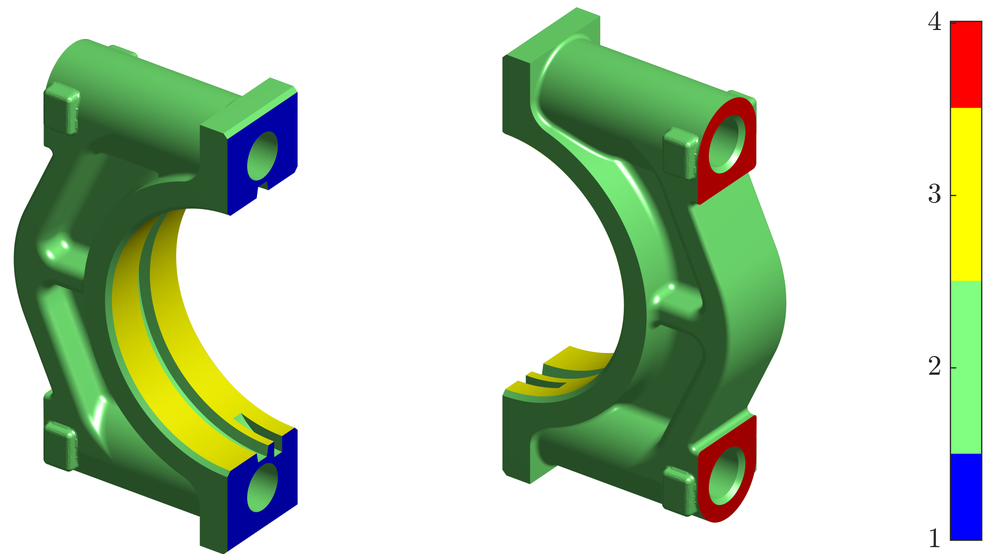}	
	\caption{Bearing cap: geometric model with the colours representing the different boundary conditions.}
	\label{fig:bearingBC}
\end{figure}

The last example considers the application of the proposed FCFV approach in a realistic setting involving the stress analysis of a bearing cap used in the automotive industry. Figure~\ref{fig:bearingBC} shows the geometry of the component, where the different colours represent the different boundary conditions. 
A homogeneous Dirichlet boundary condition is applied to the surfaces in blue, where the bearing cap is fixed. Neumann boundary conditions, enforcing a prescribed pressure of $P=130$N/mm$^2$, are imposed on the surfaces in red and yellow, corresponding to the pressure exerted by the screws and the crankshaft respectively. Homogeneous Neumann boundary conditions are imposed on the rest of the boundary surfaces in green. The bearing cap is made of cast iron with $E=130$GPa and $\nu=0.25$. 

\begin{figure}[!tb]
	\centering
	\subfigure[$u_1$]{\includegraphics[width=0.32\textwidth]{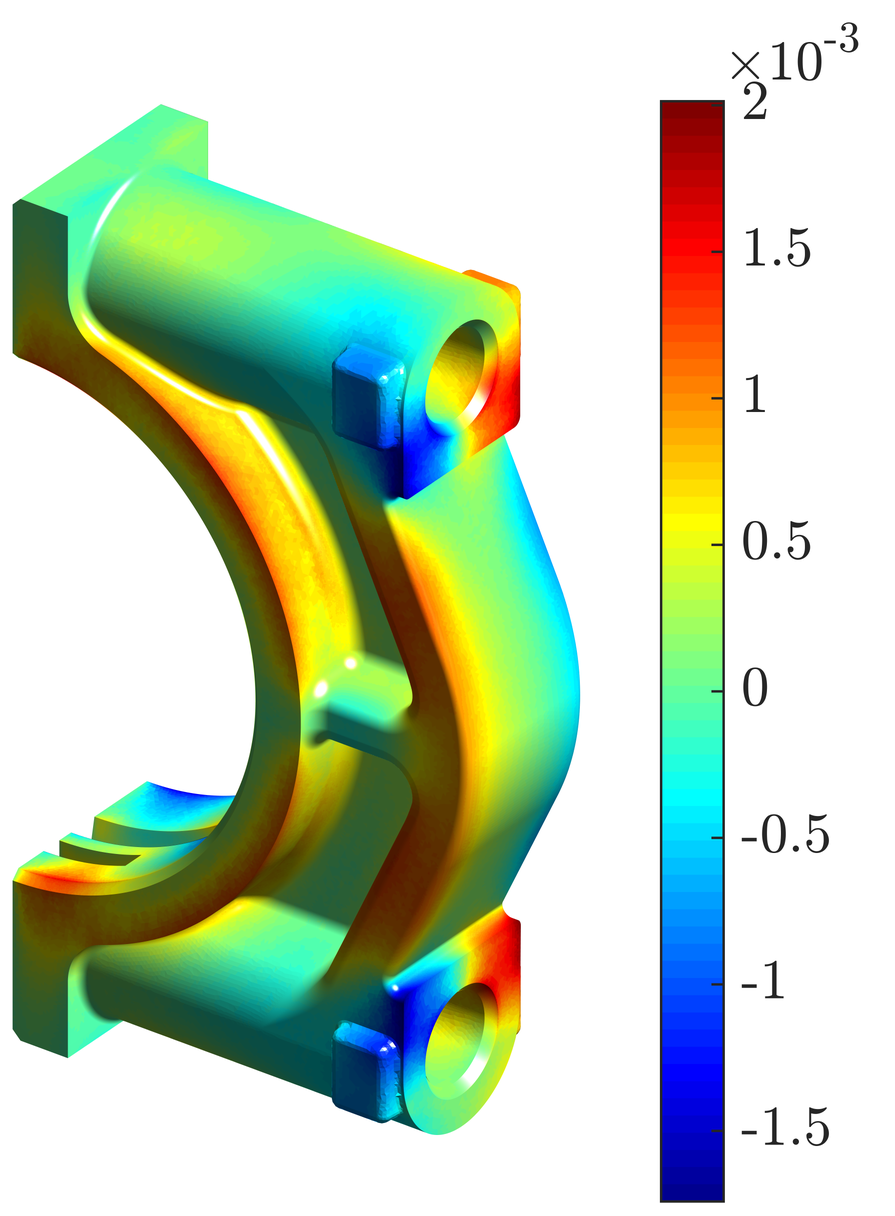}}
	\subfigure[$u_2$]{\includegraphics[width=0.32\textwidth]{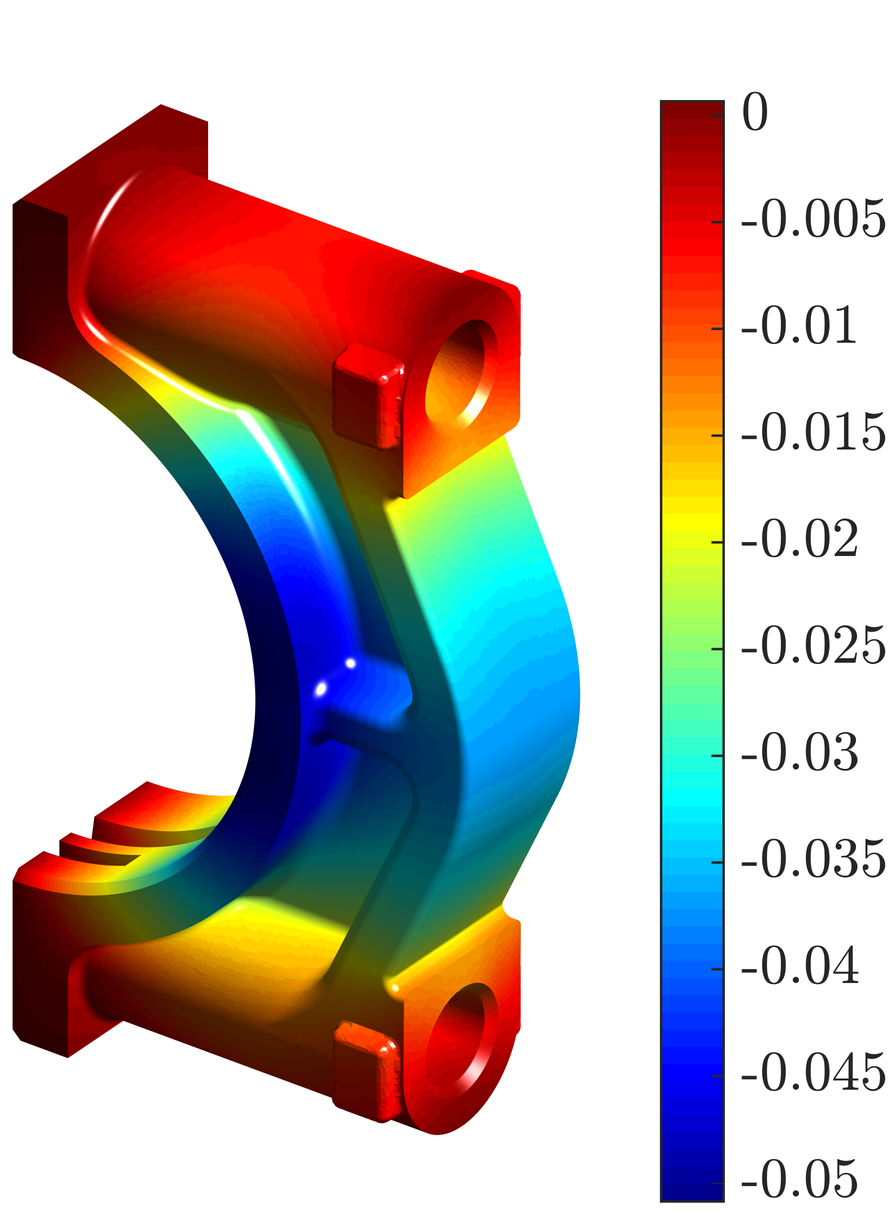}}	
	\subfigure[$u_3$]{\includegraphics[width=0.32\textwidth]{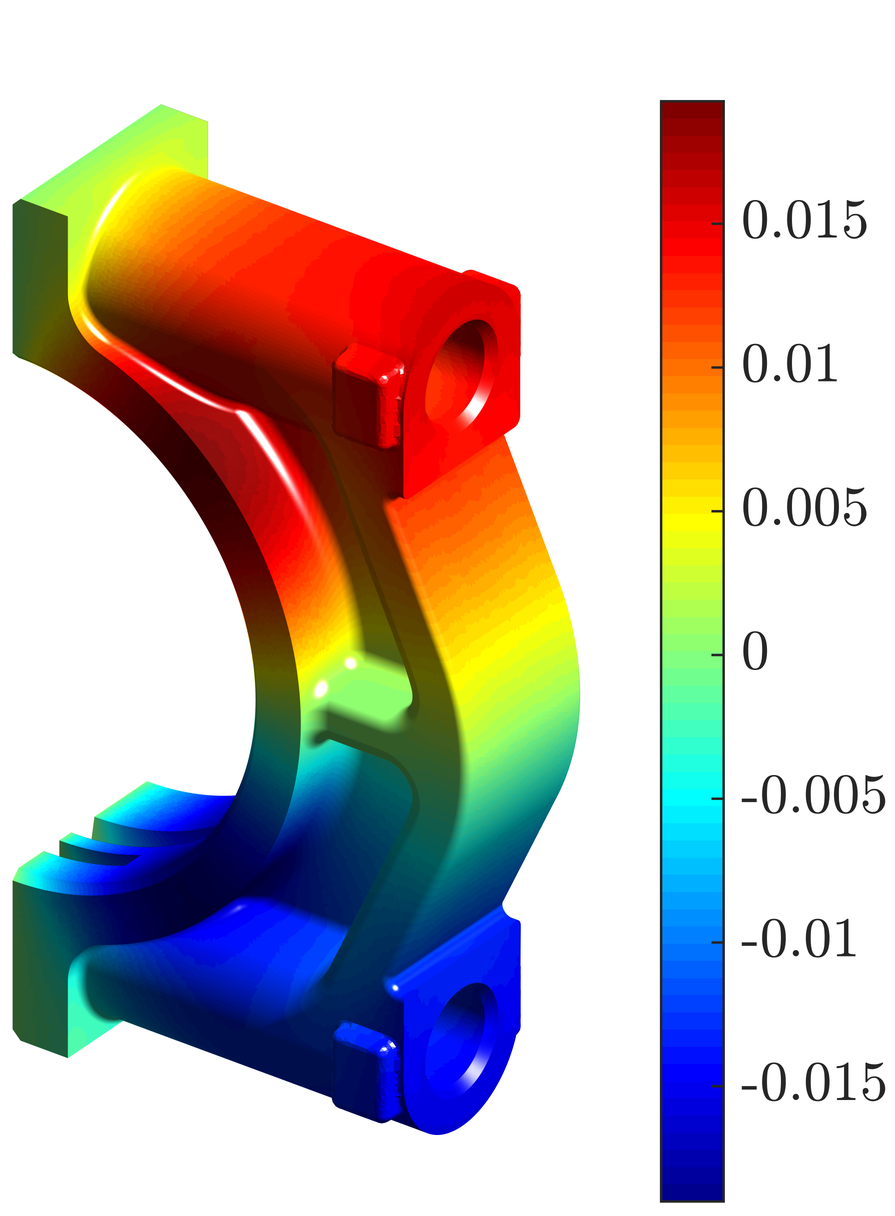}}	
	\caption{Bearing cap: components of the displacement field in mm.}
	\label{fig:bearingDisplacement}
\end{figure}

\begin{figure}[!tb]
	\centering
	\subfigure[$\sigma_{11}$]{\includegraphics[width=0.32\textwidth]{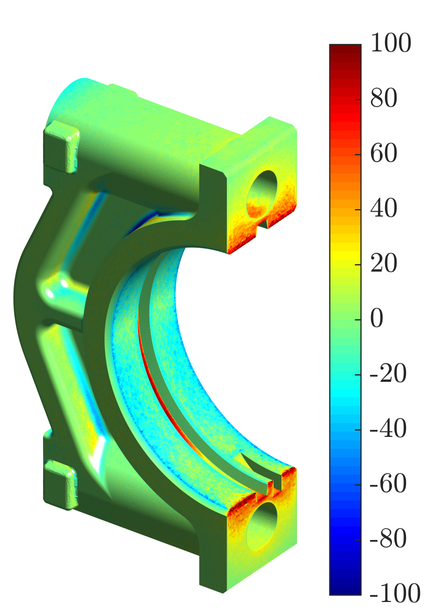}}
	\subfigure[$\sigma_{22}$]{\includegraphics[width=0.32\textwidth]{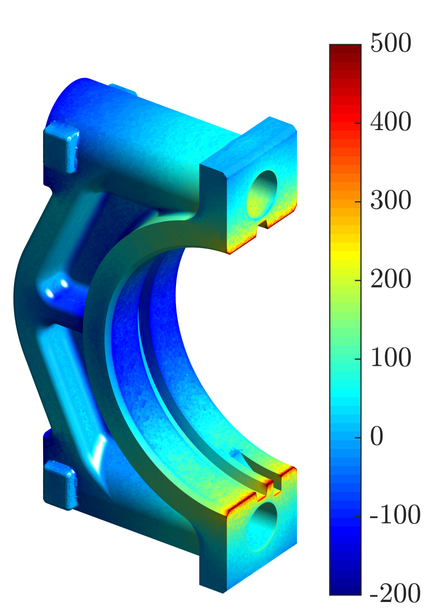}}	
	\subfigure[$\sigma_{33}$]{\includegraphics[width=0.32\textwidth]{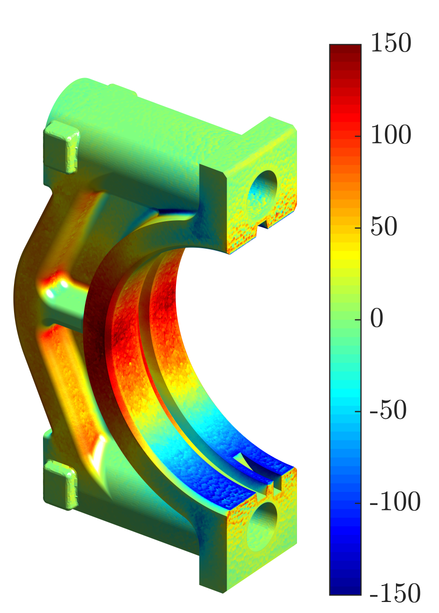}}	
	\caption{Bearing cap: three components of the stress field in N/mm$^2$.}
	\label{fig:bearingStress}
\end{figure}

The three components of the computed displacement field and the axial components of the computed stress field are represented in Figures~\ref{fig:bearingDisplacement} and \ref{fig:bearingStress} respectively. 
The computation has been performed on an unstructured mesh with 2,398,627 tetrahedral elements. The mesh has 4,256,488 internal faces and 1,081,532 external faces, leading to a global system with 15,854,985 degrees of freedom. The computation with the proposed FCFV method required 50 seconds to compute all elemental matrices, 23 seconds to perform the assembly of the global system and 41 hours to solve the system using a conjugate gradient method with no pre-conditioner. The developed code is written in Matlab and the computation was performed in an Intel$^{\tiny{\textregistered}}$
Xeon$^{\tiny{\textregistered}}$ CPU $@$ 3.70GHz and 32GB main memory available.

%==========================================================================
\section{Concluding remarks}
\label{sc:Conclusion}
%==========================================================================

A new finite volume paradigm, based on the hybridisable discontinuous Galerkin (HDG) method with constant degree of approximation, has been presented for the solution of linear elastic problems. Similar to other HDG methods, the proposed face-centred finite volume (FCFV) method provides a volumetric locking-free approach. Contrary to other HDG methods, the symmetry of the stress tensor is strongly enforced using the Voigt notation, leading to optimal convergence of the stress field components. 

The proposed FCFV method  defines the displacement unknowns on the faces (edges in two dimensions) of the mesh elements. The displacement and stress fields on each element are then recovered using closed form expressions, leading to an efficient methodology that does not require a reconstruction of the gradient of the displacement and, therefore, it is insensitive to mesh distortion. 

Numerical examples in two and three dimensions have been used to demonstrate the optimal convergence of the proposed method, its robustness when distorted meshes are considered and the absence of locking in the incompressible limit. The examples include classical benchmark test cases as well as a realistic application in three dimensions.

%==========================================================================
\section*{Acknowledgements}
%==========================================================================

This work is partially supported by the European Union's Horizon 2020 research and innovation programme under the Marie Sk\l odowska-Curie Actions (Grant number: 675919) and the Spanish Ministry of Economy and Competitiveness (Grant number: DPI2017-85139-C2-2-R).  The first author also gratefully acknowledges the financial support provided by the S\^{e}r Cymru National Research Network for Advanced Engineering and Materials (Grant number: NRN045). The second and third authors are also grateful for the financial support provided by the Generalitat de Catalunya (Grant number: 2017-SGR-1278).

%==========================================================================
%\bibliographystyle{ieeetr}
\bibliographystyle{plain}
\bibliography{Ref-HDG}
%==========================================================================

\end{document}